\documentclass[12pt]{article}
\usepackage{tikz}
\usetikzlibrary{calc,through,backgrounds}
\usepackage{multicol}

\usepackage{bm,amssymb,amsthm,amsfonts,amsmath,latexsym,cite,color, psfrag,graphicx,ifpdf,pgfkeys,pgfopts,xcolor,extarrows}

\usepackage[dvipdfm]{geometry}
\newtheorem{theo}{Theorem}[section]
\newtheorem{defi}[theo]{Definition}

\newtheorem{lem} [theo]{Lemma}
\newtheorem{coro}[theo]{Corollary}
\newtheorem{prop}[theo]{Proposition}

\newtheorem{conj}[theo]{Conjecture}

\allowdisplaybreaks

\makeatletter \@addtoreset{equation}{section}
\@addtoreset{theo}{}\makeatother

\usepackage{hyperref}

\def\qed{\hfill \rule{4pt}{7pt}}
\def\pf{\noindent {\it Proof.} }

\def\c{  \mathbf{c}}
\def\P{   \mathcal{P}  }
\def\Newton{  \mathrm{Newton}}
\begin{document}

\begin{center}

{ \bf ON THE EHRHART POLYNOMIAL OF SCHUBERT MATROIDS}

\vskip 4mm
{\footnotesize NEIL J.Y. FAN AND   YAO  LI }

\end{center}

\noindent{A{\scriptsize BSTRACT}.}
In this paper, we give a formula for the number of  lattice points in the dilations of Schubert matroid polytopes. As applications, we obtain the Ehrhart polynomials of uniform and  minimal matroids  as  special cases, and give a recursive formula for the Ehrhart polynomials of $(a,b)$-Catalan matroids.  Ferroni showed that uniform and minimal matroids are Ehrhart positive. We  show that all sparse paving Schubert matroids   are Ehrhart positive and their Ehrhart polynomials are coefficient-wisely bounded by those of  minimal and uniform matroids. This confirms a conjecture of Ferroni for the case of sparse paving Schubert matroids. Furthermore, we introduce notched rectangle matroids, which include minimal matroids, sparse paving Schubert matroids and panhandle matroids. We  show that three subfamilies of notched rectangle matroids are Ehrhart positive, and conjecture that all notched rectangle matroids are Ehrhart positive.

\section{Introduction}

Let $S$ be a subset   of $[n]:=\{1,2,\ldots,n\}$. The {\it Schubert matroid} ${\rm SM}_n(S)$ is the matroid with ground set $[n]$ and  bases
\[\{T\subseteq [n]\colon T\leq S\},\]
where $T\leq S$ means that: $|T|=|S|$ and
 the $i$-th smallest element of $T$ does not exceed that of $S$ for $1\leq i\leq |T|$. Schubert matroids were first studied by Crapo \cite{Crapo} under the name of nested matroids and rediscovered in various contexts. They have been called shifted matroids \cite{Ard-1},  freedom matroids \cite{Crapo2}, generalized Catalan matroids \cite{Bon2} and PI-matroids \cite{Billera}  in the literature.
In particular, uniform matroids, minimal matroids \cite{Fer2},  $(a,b)$-Catalan  matroids \cite{Bon} and panhandle matroids \cite{panhandle} are subclasses of Schubert matroids.
It is also worth mentioning that  Schubert matroids
are subfamilies of lattice path matroids
\cite{Bon,Bon2, Knauer}, or more generally transversal matroids \cite{Ard-1}
and positroids \cite{Oh}.

It follows from
Derksen and Fink \cite[Theorem 5.4]{Derksen}  that essentially Schubert matroids form a basis of the indicator function space of all matroids.   In a more elementary language, the Ehrhart polynomial of an arbitrary matroid polytope is an integer linear combination of Ehrhart polynomials of Schubert matroid polytopes. Moreover, Schubert matroid polytopes   are the  Minkowski summands of the Newton polytopes of key polynomials  and Schubert polynomials, or more generally, Schubitopes, see Fink, M\'esz\'aros  and St.$\,$Dizier \cite{Fin}.

Suppose that $M=([n], \mathcal{B})$ is a matroid with ground set $[n]$  and base set $\mathcal{B}$.
The {\it matroid polytope} $\P(M)$ associated to $M$ is  the convex hull
\begin{align}\label{pm}
\P(M)=\mathrm{conv}\{e_B\colon B\in \mathcal{B}\},
\end{align}
where $e_B=e_{b_1}+\cdots+e_{b_k}$ for $B=\{b_1,\ldots,b_k\}\subseteq[n]$ and $\{e_i\colon 1\leq i\leq n\}$ is the standard basis of $\mathbb{R}^n$.
Given a polytope $\P$ and a positive integer $t$, the $t$-dilation $t\P$ of $\P$ is defined as $t\P=\{t\alpha|\alpha\in \P\}$. Let $i(\P,t)=|t\P\cap \mathbb{Z}^n|$ denote the number of lattice points in $t\P$. It is well known that for integral polytopes, $i(\P, t)$ is a polynomial in $t$, called the {\it Ehrhart polynomial} of $\P$. For simplicity,  write $i(M, t)$ for $i(\P(M), t)$.

It was conjectured by De Loera, Haws and K\"{o}ppe \cite{Lorea} that all matroids are \textit{Ehrhart positive}, i.e., the Ehrhart polynomial of any matroid polytope has positive coefficients. Moreover, since matroid polytopes are specific families of generalized permutohedra, Castillo and  Liu \cite{Liu} further conjectured that generalized permutohedra are also Ehrhart positive. For the study of Ehrhart positivity of various polytopes, see the survey of Liu \cite{Liu}.
Recently, Ferroni  \cite{Fer, Fer2} showed that hypersimplices and minimal matroids are Ehrhart positive. In \cite{Fer3}, Ferroni showed that all sparse paving matroids of rank 2 are Ehrhart positive,  but provided counterexamples to both aforementioned conjectures  of all ranks greater than or equal to 3. Ferroni,   Jochemko and Schr\"{o}ter \cite{Fer4} further showed that  all   matroids of   rank 2 are  Ehrhart positive and    are coefficient-wisely bounded by  minimal and
  uniform matroids.

In this paper, we consider the Ehrhart polynomials of Schubert matroid polytopes.     We provide a formula for the number of lattice points in the $t$-dilation $t\P({\rm SM}_n(S))$ of $\P({\rm SM}_n(S))$. To this end, we  first show that $t\P({\rm SM}_n(S))$ is in fact the Newton polytope of the key polynomial $\kappa_{t\alpha}(x)$, where $\alpha$ is the indicator vector of $S$. It follows from Fink,  M\'esz\'aros  and St.$\,$Dizier \cite{Fin} that each lattice point in the Newton polytope of $\kappa_{t\alpha}(x)$ is an exponent vector of $\kappa_{t\alpha}(x)$. Then we use Kohnert algorithm to generate all the different monomials of $\kappa_{t\alpha}(x)$ and thus obtain a formula for the number of lattice points in  $t\P({\rm SM}_n(S))$.

As applications, we obtain the Ehrhart polynomials of hypersimplices  \cite{Kat} and minimal matroids \cite{Fer2} as simple special cases, and give  a recursive  formula for the Ehrhart polynomials  of  $(a,b)$-Catalan matroids. We also show that all sparse paving Schubert matroids are Ehrhart positive by proving that they are coefficient-wisely bounded by the minimal and uniform matroids.  Ferroni \cite{Fer2} conjectured that all matroids are coefficient-wisely bounded by the minimal and uniform matroids, which was disproved by Ferroni \cite{Fer3} later on. We confirm this conjecture for the case of sparse paving Schubert matroids.  Moreover, we introduce notched rectangle matroids, and show that three subfamilies of notched rectangle matroids are
Ehrhart positive. We conjecture that all notched rectangle matroids are Ehrhart positive.

To describe our results, we need some notations.
Assume that  $S\subseteq[n]$ is a finite set of positive integers. Since we only consider Schubert matroids ${\rm SM}_n(S)$, it suffices to let $n$ be the maximal element of $S$. The {\it indicator vector} $\mathbb{I}(S)$ of $S$ is the 0-1 vector $\mathbb{I}(S)=(i_1,\ldots,i_n)$, where $i_j=1$ if $j\in S$, and 0 otherwise. Clearly, $i_n=1$.  For simplicity, write $\mathbb{I}(S)=(0^{r_1},1^{r_2},\ldots,0^{r_{2m-1}},1^{r_{2m}})$,
where $0^{r_1}$ represents $r_1$ copies of 0's, $1^{r_2}$ represents $r_2$ copies of 1's, etc. Thus $S$ can be written   as an  integer sequence $r(S)=(r_1,r_2,\ldots,r_{2m})$ of length $2m$,
where $r_1\ge0$ and $r_i>0$ for $i\ge 2$.
It is easy to see that given such an integer sequence $r$, there is a unique set $S$  whose indicator vector $\mathbb{I}(S)$ can be written   in this way. We will use $S$, $r$ or $r(S)$ interchangeably  with no further clarification. For example, let $S=\{2,6,7,10\}\subseteq[10]$, then $\mathbb{I}(S)=(0,1,0^3,1^2,0^2,1)$ and $r(S)=(1,1,3,2,2,1)$.

Given $r=(r_1,r_2,\ldots,r_{2m})$, define two integer sequences $u=(u_1,\ldots,u_m)$ and $v=(v_1,\ldots,v_m)$ as follows. For $1\le i\le m$, let
\begin{align}\label{uv}
u_{i}&=\min\left\{r_{2i-1},\sum_{j=i+1}^mr_{2j}\right\} \  \ \text{and} \ \ v_{i}=\min\left\{r_{2i},\sum_{j=1}^{i-1}r_{2j-1}\right\},
\end{align}
where empty sums are interpreted as 0.
Assume that $a,b,t\ge0$ and $c\in\mathbb{Z}$ are all integers, define
\begin{align}\label{Fabct}
F(a,b,c,t)=\sum_{j=0}^{a+b}(-1)^j{a+b\choose j}{(t+1)(b-j)+a+c-1\choose a+b-1}.
\end{align}
By convention, ${0\choose0}=1$ and ${n\choose k}=0$ if $k<0$ or $n<k$. Notice that if $j>\frac{bt+c}{t+1}$ in \eqref{Fabct}, then $(t+1)(b-j)+a+c-1<a+b-1$, and thus ${(t+1)(b-j)+a+c-1\choose a+b-1}=0$.

\begin{theo}\label{main}
Let $S\subseteq [n]$ with $r(S)=(r_1,\ldots,r_{2m})$. We have
\begin{align}\label{ipt}
i({\rm SM}_n(S),t)=\sum_{(c_1,\ldots,c_m)}\prod_{j=1}^m F(r_{2j-1},r_{2j},c_{j},t),
\end{align}
where $c_1+\cdots+c_m=0$, and for $1\le j\le m$,
\begin{align*}
-tv_j\le c_j\le tu_j\ \ \text{and}\ \
c_1+\cdots+c_j\ge0.
\end{align*}
\end{theo}

Since the variable $t$ appears as the upper limit of the sum,
\eqref{ipt} is not a legitimate polynomial. Nevertheless, there are still many applications.
For example,  let $S=\{n-k+1,\ldots,n\}$,  where $n> k\ge1$. Then we obtain the uniform matroid $U_{k,n}$. In this case, $r=(n-k,k)$, $m=1, \ c_1=0$,   by Theorem \ref{main},
\[i(U_{k,n},t)=F(n-k,k,0,t).\]
The Ehrhart polynomial $i(U_{k,n},t)$  was first obtained by Katzman \cite{Kat} and then shown to have positive coefficients by Ferroni \cite{Fer}.

\begin{coro}[Katzman \cite{Kat}]\label{coro1}
We have
\[
i(U_{k,n}, t)=F(n-k,k,0,t)=\sum_{i=0}^{k-1}(-1)^i{n\choose i}{(k-i)t-i+n-1\choose n-1}.
\]
\end{coro}

Let $S=\{2,3,\ldots,k,n\}$, where $n>k\ge2$, we are lead to the minimal matroid $T_{k,n}$.
Minimal matroids were first studied independently by  Dinolt \cite{Din} and Murty \cite{Mur}. Ferroni \cite{Fer2} showed that $T_{k,n}$ is the graphic matroid of  a $(k+1)$-cycle    with one edge   replaced by $n-k$ parallel copies. In this case, $r(S)=(1,k-1,n-k-1,1)$, $u=(1,0), v=(0,1)$ and $(c_1,c_2)=(j,-j)$ for $0\le j\le t$, thus by Theorem \ref{main},
\[i(T_{k,n},t)=\sum_{j=0}^tF(1,k-1,j,t)F(n-k-1,1,-j,t).\]
Since both $F(1,k-1,j,t)$ and $F(n-k-1,1,-j,t)$ are binomials, we can re-obtain the following closed formula of $i(T_{k,n},t)$.

\begin{coro}[Ferroni \cite{Fer2}]\label{coro2}
We have
\begin{align}\label{ehr}
i(T_{k,n},t)=\frac{1}{{n-1\choose k-1}}{t+n-k\choose n-k}\sum_{j=0}^{k-1}{n-k+j-1\choose j}{t+j\choose j}.
\end{align}
\end{coro}

It is apparent that both $i(T_{k,n},t)$ and $i(T_{k,n},t-1)$ have positive coefficients.
We proceed to consider some further applications of Theorem \ref{main}.

\subsection{Notched rectangle matroids}

\begin{figure}[h]
\begin{center}
\begin{tikzpicture}

\draw [very thick](0mm,0mm)--(24mm,0mm)--(24mm,16mm)
--(48mm,16mm)--(48mm,40mm)--(0mm,40mm)--(0mm,0mm);

\draw [step=8mm] (0mm,0mm) grid (24mm,16mm);
\draw [step=8mm] (0mm,16mm) grid (48mm,40mm);
\draw [dashed,line width=2pt] (0mm,0mm)--(8mm,0mm)
--(8mm,8mm)--(16mm,8mm)--(16mm,32mm)--(40mm,32mm)
--(40mm,40mm)--(48mm,40mm);

\node at (4mm,-2mm) {$1$};\node at (12mm,-2mm) {$2$};\node at (20mm,-2mm) {$3$};

\node at (25mm,4mm) {$4$};\node at (25mm,12mm) {$5$};

\node at (10mm,4mm) {$2$};\node at (18mm,12mm) {$4$};

\node at (28mm,14mm) {$6$};\node at (36mm,14mm) {$7$};\node at (44mm,14mm) {$8$};

\node at (49mm,20mm) {$9$};\node at (50mm,28mm) {$10$};\node at (50mm,36mm) {$11$};

\node at (18mm,20mm) {$5$};\node at (18mm,28mm) {$6$};
\node at (42mm,36mm) {$10$};

\end{tikzpicture}
\end{center}
\vspace{-4mm}
\caption{A notched rectangle matroid with parameters $(3,2,3,3)$.}\label{notch}
\end{figure}
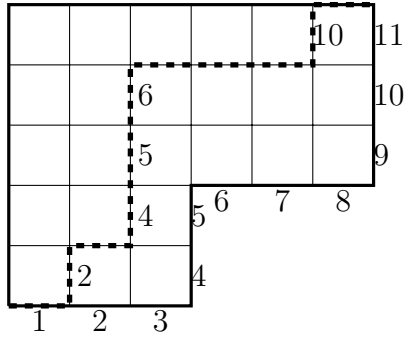

Given a lattice path $P$ from $(0,0)$ to $(m,n)$ consisting of East steps $E=(1,0)$ and North steps $N=(0,1)$, label the steps of $P$ by $1,2,\ldots,m+n$. Let $B(P)$ denote the set of labels of North steps of $P$.  For example, in Figure \ref{notch}, the dashed lattice path $P$ has $B(P)=\{2,4,5,6,10\}$. Let $U,L$ be two lattice paths from $(0,0)$ to $(m,n)$, such that $L$ never goes above $U$. The \textit{lattice path matroid} $M[U,L]$ is the matroid on the ground set $[m+n]$ with base consisting of $B(P)$, where $P$ is a lattice path from $(0, 0)$ to $(m,n)$ never going below $L$ and never going above $U$.
For the study of lattice path matroids, see    \cite{Bon,Bon2,Knauer}.

It is easy to see that a Schubert matroid ${\rm SM}_n(S)$ is a lattice path matroid $M(U,L)$, where
\[
B(U)=\{1,2,\ldots,|S|\}\ \text{and}\ B(L)=S.
\]
Bonin and   de Mier \cite[Definition 8.1]{Bon2} introduced  notch matroids, which are lattice path matroids of the form $M[U,E^mN^n]$ or $M[U,E^{m-1}NEN^{n-1}]$.
Recently,  Hanely et al. \cite{panhandle} studied another subfamily of lattice path matroids they called panhandle matroids, and conjectured that panhandle matroids are Ehrhart positive. Notice that notch matroids are not Schubert matroids in general, but panhandle matroids are in fact  Schubert matroids ${\rm SM}_n(S)$ with $r(S)=(a,b,c,1)$.

We introduce a more general family of matroids, which includes minimal matroids, sparse paving Schubert matroids and panhandle matroids.

\begin{defi}
A notched rectangle matroid is  a Schubert matroid ${\rm SM}_n(S)$ with $r(S)=(a,b,c,d)$, where $a,b,c,d$ are positive integers.
\end{defi}

Figure \ref{notch} is an illustration of the notched rectangle matroid  with parameters $(3,2,3,3)$, which is exactly the Schubert matroid ${\rm SM}_n(S)$ with $r(S)=(3,2,3,3)$, or equivalently, $S=\{4,5, 9,10,11\}$.

We express the Ehrhart polynomials of three subfamilies of notched rectangle matroids as positive combinations of $i(U_{b,a+b},t)$, which imply Ehrhart positivity of these matroids.
For convenience, let $i(r(S),t)$    denote $i({\rm SM}_n(S),t)$.

\begin{theo}  \label{spcia}
Let $a,b$ be positive integers.  Then
\begin{align}
i((a,b,a,b),t)&=\frac{1}{2}i(U_{2b,2a+2b},t)
+\frac{1}{2}i(U_{b,a+b},t)^2,
\label{ababs}\\
i((a,a,b,b),t)&=\frac{1}{2}i(U_{a+b,2a+2b},t)
+\frac{1}{2}i(U_{a,2a},t)i(U_{b,2b},t),\label{aabb}\\
i((1,1,a,a+1),t)&=i((a+1,a,1,1),t)=\frac{1}{2}(t+2)i(U_{a+1,2a+2},t).\label{11aa1}
\end{align}
\end{theo}

\begin{conj}
All notched rectangle matroids are Ehrhart positive.
\end{conj}

\subsection{$(a,b)$-Catalan matroids}

Let $r=(\overbrace{a,b,a,b,\ldots,a,b}^{2n})$, where $a,b,n\ge1$, we obtain the $(a,b)$-Catalan matroid $\mathbf{C}_n^{a,b}$ introduced by Bonin,  de Mier and  Noy \cite[Definition 3.7]{Bon}.  In particular, when $a=b=1$, we obtain the Schubert matroid $\mathbf{C}_n^{1,1}$, which is equivalent to $\text{SM}_{2n}(S)$ with
\[S=\{2,4,\ldots,2n\}.\]
In \cite{Ard-1}, Ardila studied  the Catalan matroid $\mathbf{C}_n$, which is the Schubert matroid $\text{SM}_{2n-1}(S)$ with $S=\{1,3,\ldots,2n-1\}$ and an additional loop $2n$. It is easy to see that  $\mathbf{C}_n^{1,1}$ is isomorphic to $\mathbf{C}_{n+1}$.

A composition $\sigma=(\sigma_1,\ldots,\sigma_s)$ of $n$  is an ordered nonnegative integer sequence such that $\sigma_1+\cdots+\sigma_s=n$. Let $\ell(\sigma)=s$ denote the number of parts of $\sigma$.
 Given two compositions $\sigma$ and $\sigma'$, we say that $\sigma$ and $\sigma'$ are \textit{equivalent}, denoted as $\sigma\sim\sigma'$, if $\sigma'$ can be obtained from $\sigma$ by  \textit{cyclic shifting}, i.e., $\sigma'=(\sigma_j,\ldots,\sigma_s,\sigma_{1},\ldots,\sigma_{j-1})$ for some $2\le j\le s$. Let $d(\sigma)$ denote the cardinality of the equivalence class of $\sigma$.
Denote $\Gamma_n$ by a transversal of the equivalence classes consisting of compositions of $n$ with at least two parts and minimal parts larger than 1. That is, if $\sigma\in \Gamma_n$, then $\min\{\sigma_1,\ldots,\sigma_s\}>1$, $\ell(\sigma)>1$, and if $\sigma,\sigma'\in \Gamma_n$, then $\sigma'$ and $\sigma$ are not equivalent.

\begin{theo}\label{coro3}
For $a,b\ge1$ and $n\ge2$, we have
\begin{align}\label{freef}
i(\mathbf{C}_{n}^{a,b},t)&=\frac{1}{n}i(U_{nb,na+nb},t)
-\frac{1}{n}i(U_{b,a+b},t)^n+i(U_{b,a+b},t)\cdot i(\mathbf{C}_{n-1}^{a,b},t)\nonumber\\
&\quad+\sum_{\sigma \in \Gamma _n} (-1)^{\ell(\sigma )}\frac{d(\sigma )}{\ell(\sigma )}\cdot i(\overline{\mathbf{C}}_{\sigma}^{a,b},t),
\end{align}
where $\displaystyle i(\overline{\mathbf{C}}_{\sigma}^{a,b},t)=\prod_{j=1}^{\ell(\sigma )} i(\overline{\mathbf{C}}_{\sigma_j}^{a,b},t)$ and
\begin{align}\label{cbar}
i(\overline{\mathbf{C}}_{\sigma_j}^{a,b},t)=i(\mathbf{C}_{\sigma_j}^{a,b},t)
-i(U_{b,a+b},t)\cdot i(\mathbf{C}_{\sigma_j-1}^{a,b},t),
\end{align}
and $i(\mathbf{C}_{1}^{a,b},t)=i(U_{b,a+b},t)$.
\end{theo}

For example, since $\Gamma_2=\Gamma_3=\emptyset$ and $\Gamma_4=\{(2,2)\}$, we have
\begin{align*}
i(\mathbf{C}_2^{a,b},t)&=\frac{1}{2}i(U_{2b,2a+2b},t)+\frac{1}{2}i(U_{b,a+b},t)^2\\
i(\mathbf{C}_3^{a,b},t)&=\frac{1}{3}i(U_{3b,3a+3b},t)-\frac{1}{3}i(U_{b,a+b},t)^3
+i(U_{b,a+b},t)\cdot i(\mathbf{C}_{2}^{a,b},t)\\
i(\mathbf{C}_4^{a,b},t)&=\frac{1}{4}i(U_{4b,4a+4b},t)-\frac{1}{4}i(U_{b,a+b},t)^4
+i(U_{b,a+b},t)\cdot i(\mathbf{C}_{3}^{a,b},t)+\frac{1}{2} i(\overline{\mathbf{C}}_{2}^{a,b},t).
\end{align*}
For $n=9$,  let
$
\Gamma_9=\{(7,2),(6,3),(5,4),(5,2,2),
(4,3,2),(4,2,3),(3,3,3),(3,2,2,2)\}.
$
Thus
\begin{align*}
i(\mathbf{C}_9^{a,b},t)&=\frac{1}{9}i(U_{9b,9a+9b},t)
-\frac{1}{9}i(U_{b,a+b},t)^9
+i(U_{b,a+b},t)\cdot i(\mathbf{C}_{8}^{a,b},t)\\
&\quad+i(\overline{\mathbf{C}}_{(7,2)}^{a,b},t)
+i(\overline{\mathbf{C}}_{(6,3)}^{a,b},t)
+i(\overline{\mathbf{C}}_{(5,4)}^{a,b},t)-i(\overline{\mathbf{C}}_{(5,2,2)}^{a,b},t)\\
&\quad-
i(\overline{\mathbf{C}}_{(4,3,2)}^{a,b},t)-i(\overline{\mathbf{C}}_{(4,2,3)}^{a,b},t)
-\frac{1}{3}i(\overline{\mathbf{C}}_{(3,3,3)}^{a,b},t)
+i(\overline{\mathbf{C}}_{(3,2,2,2)}^{a,b},t),
\end{align*}
where $i(\overline{\mathbf{C}}_{(4,3,2)}^{a,b},t)=i(\overline{\mathbf{C}}_{(4,2,3)}^{a,b},t).$

Computational experiments  suggest the following two conjectures.

\begin{conj}
For   integers $a,b,n\ge1$, $i(U_{nb,na+nb},t)-i(U_{b,a+b},t)^n$ has positive coefficients.
\end{conj}

\begin{conj}
For   integers $a,b,n\ge1$,
 $i(\overline{\mathbf{C}}_{n}^{a,b},t)$   has positive coefficients.
\end{conj}

Notice that the positivity of $i(\overline{\mathbf{C}}_{n}^{a,b},t)$ implies the positivity of $i(\mathbf{C}_{n}^{a,b},t)$.

\subsection{Sparse paving Schubert matroids}

Let $r=(k-1,1,1,n-k-1)$, where $n>k\ge2$, we obtain a special Schubert matroid, denoted as ${\rm Sp}_{k,n}$. In fact, as will be shown in Proposition \ref{sps},   ${\rm Sp}_{k,n}$ is a sparse paving matroid, and a Schubert matroid $\text{SM}_{n}(S)$ is sparse paving if and only if $r(S)=(n-k,k)$ or $r(S)=(k-1,1,1,n-k-1)$, namely, $\text{SM}_{n}(S)$ is a uniform matroid or
\[S=\{k,k+2,\ldots,n\}.\]

\begin{theo}\label{spp}
Sparse paving Schubert matroids are Ehrhart positive and are coefficient-wisely bounded by   minimal and uniform matroids. That is, we have the coefficient-wise inequality
\begin{align}
i(T_{k,n},t)\le i({\rm Sp}_{k,n},t)\le i(U_{k,n},t).
\end{align}
\end{theo}

The organization of this paper is as follows. In Section 2, we recall basic  definitions  and notations of matroids and key polynomials. In Section 3, we give a proof of Theorem \ref{main}.
In Section 4, we explore some further properties of $F(a,b,c,t)$ and prove Corollary \ref{coro2} and Theorem \ref{spcia}.   Section 5 is devoted to prove Theorem \ref{coro3}. Finally, we show that  sparse paving Schubert matroids are Ehrhart positive in Section 6.

\section{Preliminaries}

A {\it matroid} is a pair $M=(E, \mathcal{I})$ consisting of a finite set $E$, called the ground set,  and a collection $\mathcal{I}$  of subsets of $E$, called \textit{independent sets}, such that:
 \begin{enumerate}
        \item  [(1)] $\emptyset \in \mathcal{I}$;
        \item  [(2)] If $J\in \mathcal{I}$ and $I\subseteq J$, then $I\in \mathcal{I}$;
        \item  [(3)]  If $I, J\in \mathcal{I}$ and $|I|<|J|$, then there exists $j\in J\setminus I$
 such that $I\cup \{j\}\in \mathcal{I}$.
\end{enumerate}
By (2), a matroid $M$ is determined by the collection $\mathcal{B}$  of maximal independent sets, called the bases of $M$. By (3),  all the bases have the same cardinality, called the rank of $M$, denoted as ${\rm rk}(M)$. So we can  write $M=(E,\mathcal{B} )$. The dual   of $M$ is the matroid $M^*=(E,\mathcal{B}^*)$, where $\mathcal{B}^*=\{E\setminus B: B\in \mathcal{B}\}$. It is easy to check that the dual of  a Schubert matroid ${\rm SM}_n(S)$ is isomorphic to ${\rm SM}_n(S')$, where $r(S')$ is the reverse of $r(S)$.

A subset $I$ of $E$ is called  \textit{dependent}  if it is not an independent set. If $C\subseteq E$ is dependent but every proper subset of $C$ is independent, we say that $C$ is a \textit{circuit}. A  subset $F$ of $E$  is called a \textit{flat} if  ${\rm rk}_M(F\cup \{a\}) > {\rm rk}_M(F)$  for every $a\notin F$. A \textit{hyperplane} $H$   is a flat such that ${\rm rk}_M(H)={\rm rk}(M)-1$.

We say that $M$ is \textit{paving} if every
circuit of $M$ has cardinality at least ${\rm rk}(M)$. A matroid $M$ is  \textit{sparse paving} if both $M$ and its dual  are paving.
A matroid   is sparse paving if and only if every subset of cardinality ${\rm rk}(M)$ is either a basis or a circuit-hyperplane, see, for example, Bonin \cite{Bon1} or Ferroni \cite[Lemma 2.7]{Fer3}.

The rank  function ${\rm rk}_M: 2^E\rightarrow \mathbb{Z}$  of $M$  is  defined by
 \[{\rm rk}_M(T)=\max\{|T\cap B|\colon B\in \mathcal{B}\},\ \ \ \text{for $T\subseteq E$}.\]
Let ${\rm rk}_S$ denote the rank function of a Schubert matroid ${\rm SM}_n(S)$.
Fan and Guo \cite[Theorem 3.3]{fg} provided an efficient algorithm to compute ${\rm rk}_S(T)$ for any $T\subseteq[n]$.
It is well known that the matroid polytope $\P(M)$ defined in \eqref{pm} associated to a matroid $M=([n],\mathcal{B})$ is a generalized permutohedron perametrized by the rank function of $M$, see, for example, Fink, M\'esz\'aros and St.$\,$Dizier \cite{Fin}.  To be specific,
\begin{equation}\label{DE}
\P(M)=\left\{x\in \mathbb{R}^n\colon \sum_{i\in [n]}  x_i={\rm rk}_M([n])
\ \  \text{and}\ \ \sum_{i\in T}x_i\leq  {\rm rk}_M(T)\ \ \text{for $T\subsetneq [n]$}\right\}.
\end{equation}

The {\it key polynomials}  $\kappa_\alpha(x)$ associated to compositions $\alpha\in \mathbb{Z}_{\geq 0}^n$
can be defined recursively as below.  If $\alpha=(\alpha_1,\alpha_2,\ldots, \alpha_n)$ is a partition (i.e., weakly decreasing), then set
$\kappa_\alpha(x)=x_1^{\alpha_1}x_2^{\alpha_2} \cdots x_n^{\alpha_n}.$
Otherwise, choose an index  $i$ such that $\alpha_i<\alpha_{i+1}$, and let $\alpha'$ be
 obtained from $\alpha$ by interchanging $\alpha_i$ and $\alpha_{i+1}$. Set
\begin{equation*}
\kappa_\alpha(x)=\partial_i(x_i\kappa_{\alpha'}(x)).
\end{equation*}
Here $\partial_i$ is the  divided difference operator sending
  a polynomial $f(x)\in \mathbb{R}[x_1,\ldots,x_n]$
  to
  \[\partial_i(f(x))=\frac{f(x)-s_i f(x)}{x_i-x_{i+1}},\]
where $s_i f(x)$ is obtained from $f(x)$ by interchanging $x_i$ and $x_{i+1}$.
Key polynomials are also called   Demazure characters, they are characters of the Demazure
modules for the general linear groups, see Demazure \cite{Dem-1,Dem-2}.

\begin{figure}[h]
\begin{center}
\begin{tikzpicture}[scale=0.5]
      \draw [-, shift={(0,0)}] (0,0)--(0,4);
      \filldraw [fill=gray, fill opacity=0.3, shift={(0,0)}] (0,0) rectangle (1,1);
      \filldraw [fill=gray, fill opacity=0.3, shift={(0,0)}] (1,0) rectangle (2,1);
      \filldraw [fill=gray, fill opacity=0.3, shift={(0,0)}] (0,2) rectangle (1,3);
      \filldraw [fill=gray, fill opacity=0.3, shift={(0,0)}] (1,2) rectangle (2,3);
      \filldraw [fill=gray, fill opacity=0.3, shift={(0,0)}] (2,2) rectangle (3,3);
      \filldraw [fill=gray, fill opacity=0.3, shift={(0,0)}] (0,3) rectangle (1,4);
\end{tikzpicture}
\end{center}
\vspace{-6mm}
\caption{ The skyline diagram $D(\alpha)$ for  $\alpha=(1,3,0,2)$.}
\label{RS}
\end{figure}
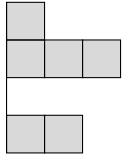

Kohnert \cite{Koh} found that the key polynomial $\kappa_\alpha(x)$  can be generated by applying  the {\it Kohnert algorithm} to the    skyline diagram  of $\alpha$, see also Reiner and Shimozono \cite{Rei}.
Recall that the {\it  skyline diagram} $D(\alpha)$ of a composition $\alpha=(\alpha_1,\ldots,\alpha_n)$  is a diagram consisting of the first
$\alpha_i$ boxes in row $i$.
For example,
Figure \ref{RS} is the skyline diagram of $\alpha=(1,3,0,2)$.

The Kohnert  algorithm is defined based on Kohnert moves on diagrams.
A {\it  diagram} $D$ is a finite collection of boxes in   $\mathbb{Z}_{>0}^2$.
A  box in row $i$ and column $j$ of the grid is denoted  $(i,j)$.
Here,  the rows  (respectively, columns) are labeled  increasingly from top to bottom (respectively, from left to right).
 A {\it  Kohnert move} on  $D$
selects the rightmost  box in a row  of $D$ and moves it within its column
 up to the first available position. To be specific,   a box $(i,j)$ of $D$ can be moved up to  a position $(i', j)$ by a
Kohnert move whenever: (i) the box $(i,j)$ is the rightmost box in the $i$-th row of $D$,
(ii) the box $(i',j)$ does not belong to $D$, and (iii) for any $i'<r<i$, the box $(r, j)$ belongs to $D$.

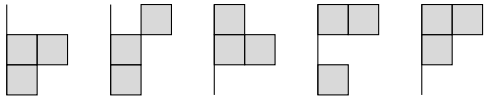
\begin{figure}[h]
\begin{center}
      \begin{tikzpicture}[scale=0.4]
      \draw (0,0) -- (0,3);
      \filldraw [fill=gray, fill opacity=0.3, shift={(0,0)}] (0,0) rectangle (1,1);
      \filldraw [fill=gray, fill opacity=0.3, shift={(0,0)}] (0,1) rectangle (1,2);
      \filldraw [fill=gray, fill opacity=0.3, shift={(0,0)}] (1,1) rectangle (2,2);
      \end{tikzpicture}
      \quad
      \begin{tikzpicture}[scale=0.4]
      \draw (0,0) -- (0,3);
      \filldraw [fill=gray, fill opacity=0.3, shift={(0,0)}] (0,0) rectangle (1,1);
      \filldraw [fill=gray, fill opacity=0.3, shift={(0,0)}] (0,1) rectangle (1,2);
      \filldraw [fill=gray, fill opacity=0.3, shift={(0,0)}] (1,2) rectangle (2,3);
      \end{tikzpicture}
      \quad
      \begin{tikzpicture}[scale=0.4]
      \draw (0,0) -- (0,3);
      \filldraw [fill=gray, fill opacity=0.3, shift={(0,0)}] (0,2) rectangle (1,3);
      \filldraw [fill=gray, fill opacity=0.3, shift={(0,0)}] (0,1) rectangle (1,2);
      \filldraw [fill=gray, fill opacity=0.3, shift={(0,0)}] (1,1) rectangle (2,2);
      \end{tikzpicture}
      \quad
      \begin{tikzpicture}[scale=0.4]
      \draw (0,0) -- (0,3);
      \filldraw [fill=gray, fill opacity=0.3, shift={(0,0)}] (0,0) rectangle (1,1);
      \filldraw [fill=gray, fill opacity=0.3, shift={(0,0)}] (0,2) rectangle (1,3);
      \filldraw [fill=gray, fill opacity=0.3, shift={(0,0)}] (1,2) rectangle (2,3);
      \end{tikzpicture}
      \quad
      \begin{tikzpicture}[scale=0.4]
      \draw (0,0) -- (0,3);
      \filldraw [fill=gray, fill opacity=0.3, shift={(0,0)}] (0,1) rectangle (1,2);
      \filldraw [fill=gray, fill opacity=0.3, shift={(0,0)}] (0,2) rectangle (1,3);
      \filldraw [fill=gray, fill opacity=0.3, shift={(0,0)}] (1,2) rectangle (2,3);
      \end{tikzpicture}

\end{center}
\vspace{-3mm}
\caption{Kohnert diagrams for $\alpha=(0,2,1)$.}\label{Koh}
\end{figure}

A  {\it  Kohnert diagram} for $D(\alpha)$ is the diagram obtained from $D(\alpha)$
by applying a sequence of Kohnert moves. For a diagram $D$, let
$x^D=\prod_{(i,j)\in D}x_i$.
Kohnert \cite{Koh} showed that
\begin{equation*}
\kappa_\alpha(x)=\sum_{D}x^D,
\end{equation*}
where  the sum takes over all the Kohnert diagrams for $D(\alpha)$. For example, Figure \ref{Koh} displays   all the Kohnert diagrams for $\alpha=(0,2,1)$. Thus $\kappa_{(0,2,1)}(x)=x_2^2x_3+x_1x_2x_3
+x_1x_2^2+x_1^2x_3+x_1^2x_2$.

\section{Proof of Theorem \ref{main}}

In order to give  a proof of Theorem \ref{main}, we first show that  $t\P({\rm SM}_n(S))$  is in fact the Newton polytope of a key polynomial.

Recall that for a polynomial
\[f=\sum_{\alpha\in \mathbb{Z}_{\geq 0}^n}c_\alpha x^\alpha \in \mathbb{R}[x_1,\ldots,x_n],\]
  the \textit{Newton polytope}  of $f$ is the convex hull of the exponent vectors of $f$:
\[\Newton(f)=\mathrm{conv}(\{\alpha\colon c_\alpha\neq 0\}).\]
It is obvious that each exponent vector of $f$ is a lattice point in $\Newton(f)$. Monical, Tokcan and Yong \cite{Mon} introduced the notion of \textit{saturated Newton polytope} (SNP) of a polynomial $f$, i.e.,  $f$ has saturated Newton polytope (SNP)  if every lattice point in $\Newton(f)$ is also an  exponent vector of $f$.
It was conjectured by Monical, Tokcan and Yong \cite{Mon} and  proved by Fink,   M\'esz\'aros  and   St.$\,$Dizier \cite{Fin} that  key polynomials $\kappa_{\alpha}(x)$ have  SNP.

Moreover, Fink, M\'esz\'aros, St.$\,$Dizier \cite{Fin} also showed that the Newton polytopes of key polynomials $\kappa_{\alpha}(x)$ are the Minkowski sum of Schubert matroid polytopes associated to the columns of $D(\alpha)$.  More precisely,   let $D(\alpha)=(D_1,\ldots, D_n)$,
where  $D_j$ is the $j$-th column of $D(\alpha)$. View  $D_j$ as  a subset of $[n]$:
\[D_j=\{1\leq i\leq n\colon (i,j)\in D_j\}.\]
Then the column $D_j$ defines a Schubert matroid ${\rm SM}_n(D_j)$. Let ${\rm rk}_j$ denote the rank function of
${\rm SM}_n(D_j)$.
Then
\begin{align}\label{ntkey}
\Newton(\kappa_{\alpha})&=\P({\rm SM}_n(D_1))+\cdots+\P({\rm SM}_n(D_n))\nonumber\\[5pt]
&=\left\{x\in \mathbb{R}^n\colon \sum_{i\in [n]}  x_i=|D(\alpha)|
\ \  \text{and}\ \ \sum_{i\in T}x_i\leq  {\rm rk}_{\alpha}(T) \ \ \text{for $T\subsetneq [n]$}\right\},
\end{align}
where $|D(\alpha)|$ denotes the number of boxes in $D(\alpha)$ and
\[{\rm rk}_{\alpha}(T)={\rm rk}_1(T)+\cdots+{\rm rk}_n(T).\]

\begin{lem} \label{guodu}
Let $S$ be a subset of $[n]$ and $\alpha=\mathbb{I}(S)$ be the indicator vector of $S$. Given any positive integer $t$, we have
\[
t\P({\rm SM}_n(S))=\Newton(\kappa_{t\alpha}).
\]
\end{lem}

\pf   It is easy to see that ${\rm rk}_{S}([n])=|S|$ is the number of elements in $S$.
By \eqref{DE}, we find  that
\[t\P({\rm SM}_n(S))=
\left\{x\in \mathbb{R}^n\colon \sum_{i\in [n]}  x_i=t\cdot |S|
\ \  \text{and}\ \ \sum_{i\in T}x_i\leq t\cdot {\rm rk}_{S}(T)\ \ \text{for $T\subsetneq [n]$}\right\}.
\]
On the other hand, since now $\alpha=\mathbb{I}(S)$ is a 0-1 vector, $D(t\alpha)$ has exactly $t$ columns, every column  determines the same Schubert matroid, which is exactly ${\rm SM}_n(S)$. Moreover, $|D(t\alpha)|=t\cdot|S|$ and ${\rm rk}_{t\alpha}(T)=t\cdot {\rm rk}_S(T)$.   Thus by \eqref{ntkey}, we conclude that
\[
\Newton(\kappa_{t\alpha})=t\P({\rm SM}_n(S)).
\]
This completes the proof. \qed

Now we are in a position to give a proof of Theorem \ref{main}.

\noindent{\it Proof of Theorem \ref{main}}.
By Lemma \ref{guodu}, the number of lattice points in $t\P({\rm SM}_n(S))$   is the same as that in $\Newton(\kappa_{t\alpha})$. Since key polynomials have saturated Newton polytopes, $i(\P({\rm SM}_n(S)),t)$ is equal to the number of different monomials in $\kappa_{t\alpha}(x)$.

Now we enumerate all the different monomials in  $\kappa_{t\alpha}(x)$ by Kohnert algorithm. Let $D(t\alpha)$ be the skyline diagram of $t\alpha$.
Let $D$ be a Kohnert diagram obtained from $D(t\alpha)$ by applying a sequence of Kohnert moves. Let $n=r_1+r_2+\cdots+r_{2m}$ denote the number of parts of $\alpha$, or equivalently, the number of rows of $D$, and  denote $\beta=(\beta_1,\ldots,\beta_n)$, where $\beta_i$ is the number of boxes in the $i$-th row of $D$. Clearly, we have $0\le\beta_i\le t$.
For $1\le j\le m$, let
\[
d_j=r_1+r_2+\cdots+r_{2j}
\]
and
\[c_j=\sum_{i=d_{j-1}+1}^{d_j}\beta_i-r_{2j}t,\]
where $d_0=0$. Since the $\beta_1+\cdots+\beta_n=(r_2+r_4+\cdots+r_{2m})t$,
we have
\[c_1+c_2+\cdots+c_m=(\beta_1+\cdots+\beta_n)-(r_2+r_4+\cdots+r_{2m})t=0.\]
It is also easy to see that the number of boxes in the top $d_j$ rows of $D$ is larger than or equal to that of $D(t\alpha)$, and the number of boxes in the  bottom $d_m-d_j$ rows of $D$ is smaller than that of $D(t\alpha)$. That is,
\[
\beta_1+\cdots+\beta_{d_j}\ge (\alpha_1+\cdots+\alpha_{d_j})t= (r_2+r_4+\cdots+r_{2j})t.
\]
and
\[
\beta_{d_j+1}+\cdots+\beta_{d_m}\le (\alpha_{d_j+1}+\cdots+\alpha_{d_m})t
=(r_{2j+2}+\cdots+r_{2m})t.
\]
Thus we have
\[
c_1+c_2+\cdots+c_j=(\beta_1+\cdots+\beta_{d_j})
-(r_2+r_4+\cdots+r_{2j})t\ge0.
\]
Moreover, we have
\begin{align*}
c_j&=\sum_{i=d_{j-1}+1}^{d_j}\beta_i-r_{2j}t\le
\sum_{i=d_{j-1}+1}^{d_m}\beta_i-r_{2j}t\\
&\le\sum_{i=d_{j-1}+1}^{d_m}t\alpha_i-r_{2j}t
=(r_{2j+2}+\cdots+r_{2m})t.
\end{align*}
And
\begin{align*}
c_j&=\sum_{i=d_{j-1}+1}^{d_j}\beta_i-r_{2j}t\le
\sum_{i=d_{j-1}+1}^{d_j}t-r_{2j}t=(r_{2j-1}+r_{2j})t-r_{2j}t=r_{2j-1}t.
\end{align*}
Thus we have
\begin{align}
c_j\le \min\{r_{2j-1},r_{2j+2}+\cdots+r_{2m}\}t=tu_j.
\end{align}

Similarly, we have
\begin{align*}
c_j&=\sum_{i=1}^{d_j}\beta_i-\sum_{i=1}^{d_{j-1}}\beta_i-r_{2j}t
\ge\sum_{i=1}^{d_j}t\alpha_i-\sum_{i=1}^{d_{j-1}}\beta_i-r_{2j}t\\[5pt]
&\ge\sum_{i=1}^{d_j}t\alpha_i-\sum_{i=1}^{d_{j-1}}t-r_{2j}t
=\left(\sum_{i=1}^{j}r_{2i}\right)t-\sum_{i=1}^{d_{j-1}}t-r_{2j}t\\[5pt]
&=-(r_1+r_3+\cdots+r_{2j-3})t.
\end{align*}
And
\[
c_j=\sum_{i=d_{j-1}+1}^{d_j}\beta_i-r_{2j}t\ge-r_{2j}t.
\]
Then we find
\begin{align}
c_j\ge-\min\{r_{2j},r_1+r_3+\cdots+r_{2j-3}\}t=-tv_j.
\end{align}
Therefore, $\beta=(\beta_1,\ldots,\beta_n)$ satisfies the following system of equations
\begin{align}\label{eqt}
\left\{
  \begin{array}{c}
  \ \  \sum\limits_{i=1}^{d_1}\beta_i =r_2t+c_1,   \\[10pt]
    \sum\limits_{i=d_1+1}^{d_2}\beta_i=r_4t+c_2,   \\[10pt]
    \vdots\\
 \sum\limits_{i=d_{m-1}+1}^{d_m}\beta_i=r_{2m}t+c_m,
  \end{array}
\right.
\end{align}
where  $c_1+c_2+\cdots+c_m=0$, and for $1\le j\le m$,
\begin{align*}
-tv_j\le c_j\le tu_j\ \ \text{and}\ \
c_1+c_2+\cdots+c_j\ge0.
\end{align*}

Now we enumerate the number of nonnegative integer solutions of the equation
\begin{align}\label{eqj}
\sum_{i=d_{j-1}+1}^{d_j}\beta_i=r_{2j}t+c_j, \quad (0\le \beta_i\le t).
\end{align}
Since $0\le \beta_i\le t$, it is easy to see that the number of solutions of   equation \eqref{eqj} is the coefficient of $x^{r_{2j}t+c_j}$ in
\begin{align}
&(1+x+\cdots+x^t)^{d_j-d_{j-1}}\nonumber\\
&=(1+x+\cdots+x^t)^{r_{2j-1}+r_{2j}}\nonumber\\
&=(1-x^{t+1})^{r_{2j-1}+r_{2j}}\cdot(1-x)^{-(r_{2j-1}+r_{2j})}\nonumber\\
&=\left(\sum_{i'=0}^{r_{2j-1}+r_{2j}}(-1)^{i'}{r_{2j-1}+r_{2j} \choose i'}x^{(t+1)i'}\right)\left(\sum_{j'=0}^\infty{r_{2j-1}+r_{2j}+j'-1\choose j'}x^{j'}\right)\nonumber\\
&=\sum_{j'=0}^\infty\sum_{i'=0}^{r_{2j-1}+r_{2j}}(-1)^{i'}{r_{2j-1}+r_{2j} \choose i'}{r_{2j-1}+r_{2j}+j'-1\choose j'}x^{j'+(t+1)i'}.\label{zj}
\end{align}
Let $j'=r_{2j}t+c_j-(t+1)i'$ in \eqref{zj}, we see that  the coefficient of $x^{r_{2j}t+c_j}$ is
\begin{align}\label{origin}
F(r_{2j-1},r_{2j},c_{j},t):=
\sum_{i'=0}^{r_{2j-1}+r_{2j}}(-1)^{i'}{r_{2j-1}+r_{2j} \choose i'}{(t+1)(r_{2j}-i')+r_{2j-1}+c_j-1\choose r_{2j-1}+r_{2j}-1}.
\end{align}
Consequently, the number of different monomials in $\kappa_{t\alpha}(x)$ is
\[
\sum_{(c_1,\ldots,c_m)}\prod_{j=1}^m F(r_{2j-1},r_{2j},c_{j},t).
\]

Conversely, suppose that $(\beta_1,\ldots,\beta_n)$ is an integer sequence such that $0\le \beta_i\le t$  and $(\beta_1,\ldots,\beta_n)$ satisfies the system of equations \eqref{eqt}, we shall show that there is a diagram $D$ whose $i$-th row has $\beta_i$ boxes and $D$  can be obtained from $D(t\alpha)$ by  applying   Kohnert moves.

\begin{figure}[h]
\begin{center}
\begin{tikzpicture}[scale=0.5]
\draw (0,0) -- (0,6);
\draw [-, shift={(0,0)}] (0,0)--(0,4);
      \filldraw [fill=gray, fill opacity=0, shift={(0,0)}] (0,0) rectangle (1,1);
      \filldraw [fill=gray, fill opacity=0, shift={(0,0)}] (1,0) rectangle (2,1);
      \filldraw [fill=gray, fill opacity=0, shift={(0,0)}] (2,0) rectangle (3,2);
      \filldraw [fill=gray, fill opacity=0, shift={(0,0)}] (0,1) rectangle (1,2);
      \filldraw [fill=gray, fill opacity=0, shift={(0,0)}] (1,1) rectangle (2,2);
      \filldraw [fill=gray, fill opacity=0, shift={(0,0)}] (2,1) rectangle (3,2);
      \filldraw [fill=gray, fill opacity=0, shift={(0,0)}] (0,3) rectangle (1,4);
      \filldraw [fill=gray, fill opacity=0, shift={(0,0)}] (1,3) rectangle (2,4);
      \filldraw [fill=gray, fill opacity=0, shift={(0,0)}] (2,3) rectangle (3,4);
\draw (5,2.5) node {$\longrightarrow$};
%\draw (1.4,-.7) node {$D(t\alpha)$};
\end{tikzpicture}
\quad
\begin{tikzpicture}[scale=0.5]
\draw (0,0) -- (0,6);
\draw [-, shift={(0,0)}] (0,0)--(0,4);
      \filldraw [fill=gray, fill opacity=0, shift={(0,0)}] (0,0) rectangle (1,1);
      \filldraw [fill=gray, fill opacity=0, shift={(0,0)}] (1,0) rectangle (2,1);
      \filldraw [fill=gray, fill opacity=0, shift={(0,0)}] (2,0) rectangle (3,2);
      \filldraw [fill=gray, fill opacity=0, shift={(0,0)}] (0,1) rectangle (1,2);
      \filldraw [fill=gray, fill opacity=0, shift={(0,0)}] (1,1) rectangle (2,2);
      \filldraw [fill=gray, fill opacity=0, shift={(0,0)}] (2,1) rectangle (3,2);
      \filldraw [fill=gray, fill opacity=0, shift={(0,0)}] (0,3) rectangle (1,4);
      \filldraw [fill=gray, fill opacity=0, shift={(0,0)}] (1,3) rectangle (2,4);
      \filldraw [fill=gray, fill opacity=0, shift={(0,0)}] (2,3) rectangle (3,4);
      \draw (0.5,3.5) node {$2$};\draw (1.5,3.5) node {$1$};\draw (2.5,3.5) node {$1$};
      \draw (0.5,1.5) node {$4$};\draw (1.5,1.5) node {$3$};\draw (2.5,1.5) node {$2$};
      \draw (0.5,.5) node {$6$};\draw (1.5,.5) node {$4$};\draw (2.5,.5) node {$4$};
\draw (5,2.5) node {$\longrightarrow$};

\end{tikzpicture}
\quad
\begin{tikzpicture}[scale=0.5]
\draw (0,0) -- (0,6);
\draw [-, shift={(0,0)}] (0,0)--(0,4);
      \filldraw [fill=gray, fill opacity=1, shift={(0,0)}] (0,0) rectangle (1,1);
      \filldraw [fill=gray, fill opacity=1, shift={(0,0)}] (1,2) rectangle (2,3);
      \filldraw [fill=gray, fill opacity=1, shift={(0,0)}] (2,2) rectangle (3,3);
      \filldraw [fill=gray, fill opacity=1, shift={(0,0)}] (0,2) rectangle (1,3);
      \filldraw [fill=gray, fill opacity=1, shift={(0,0)}] (1,3) rectangle (2,4);
      \filldraw [fill=gray, fill opacity=1, shift={(0,0)}] (2,4) rectangle (3,5);
      \filldraw [fill=gray, fill opacity=1, shift={(0,0)}] (0,4) rectangle (1,5);
      \filldraw [fill=gray, fill opacity=1, shift={(0,0)}] (1,6) rectangle (2,5);
      \filldraw [fill=gray, fill opacity=1, shift={(0,0)}] (2,6) rectangle (3,5);

\end{tikzpicture}

\end{center}
\vspace{-5mm}
\caption{An illustration of the construction of $D$. }\label{figD}
\end{figure}
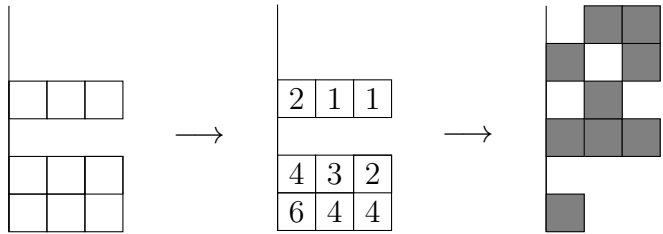

First of all, by adding all the equations in \eqref{eqt} together and combing the condition $c_1+\cdots+c_m=0$,  we have \[\beta_1+\beta_2+\cdots+\beta_n=(r_2+r_4+\cdots+r_{2m})t.\]
That is, $\beta_1+\beta_2+\cdots+\beta_n$ is equal to the total number of boxes in $D(t\alpha)$.
We construct $D$ as follows. Fill the sequence of integers $1^{\beta_1},2^{\beta_2},\ldots,n^{\beta_n}$ into the boxes of $D(t\alpha)$ along the rows from top to bottom and from right to left. Then move the box $(i,j)$ filled with $k$ to $(k,j)$. Denote the resulting diagram by $D$.  For example, Figure \ref{figD} displays the construction of $D$ for $t\alpha=(0,0,3,0,3,3)$ and $\beta=(2,2,1,3,0,1)$.

We aim to show that $D$ is indeed a Kohnert diagram.
Since $0\le\beta_i\le t$, by the construction of $D$, it is easy to see that there do not exist two boxes in the same column of $D(t\alpha)$ that are  filled with the same integer. By the definition of Kohnert moves, to show that $D$ is indeed a Kohnert diagram, it suffices to show that there does not exist a box of $D(t\alpha)$ which is filled with an integer larger than its row index.

Suppose to the contrary that there is a box $(i,j)$ filled with $s$ and $s>i$. Without loss of generality, we can  assume that $(i,j)$ is such a box with $i$ largest.
Since $i<s\le n=r_1+\cdots+r_{2m}$, there exist  integers $k,l$ such that $k\le l$ and
\begin{align}
r_1+r_2+\cdots+r_{2k-1}&<i\le r_1+r_2+\cdots+r_{2k}=d_k\\
r_1+r_2+\cdots+r_{2l-1}&<s\le r_1+r_2+\cdots+r_{2l}=d_l.
\end{align}
There are three cases.

Case 1. $l=k$, that is,
\[
r_1+r_2+\cdots+r_{2l-1}<i<s\le r_1+r_2+\cdots+r_{2l}=d_l.
\]
By \eqref{eqt}, we have
\[
\beta_1+\cdots+\beta_{d_l}
=(r_2+r_4+\cdots+r_{2l})t+c_1+\cdots+c_l.
\]
Since $c_1+\cdots+c_l\ge0$, we find that
\[
\beta_1+\cdots+\beta_{d_l}\ge (r_2+r_4+\cdots+r_{2l})t.
\]
That is to say, the integers $1,2,\ldots,d_l$ must occupy at least the top $d_l=r_1+r_2+\cdots+r_{2l}$  rows (including empty rows) of $D(t\alpha)$. On the other hand, since the box $(i,j)$ is filled with $s$ and $s>i$, we see that the rows $i+1,\ldots,s,s+1,\ldots,d_l$ of $D(t\alpha)$ are occupied by some of the integers among $s,s+1,\ldots,d_l$. In particular, the integers $s,s+1,\ldots,d_l$ must occupy at least  $d_l-s+1$ rows of $D(t\alpha)$.  Therefore,
\[
\beta_s+\beta_{s+1}+\cdots+\beta_{d_l}>(d_l-s+1)t.
\]
Thus there must exist some $\beta_j>t$, which contradicts with the assumption $0\le \beta_j\le t$.

Case 2.  $l\ge k+1$ and $i=d_k$. In this case, we have
\[s>r_1+r_2+\cdots+r_{2l-1}\ge r_1+r_2+\cdots+r_{2k+1}=i+r_{2k+1}.\]
Similar to Case 1, we see that the integers $s,s+1,\ldots,d_l$ must occupy all the boxes  of $D(t\alpha)$ in the rows $s,s+1,\ldots,d_l$. Thus
\[
\beta_s+\beta_{s+1}+\cdots+\beta_{d_l}>(d_l-s+1)t,
\]
which is a contradiction.

Case 3.  $l\ge k+1$ and $i<d_k$. In this case, we have $i+1\le d_k$ and
\[s>r_1+r_2+\cdots+r_{2l-1}\ge r_1+r_2+\cdots+r_{2k+1}\ge i+r_{2k+1}\ge i+1.\]
By the choice of $(i,j)$,  $i$ is the largest index such that $s>i$ and $(i,j)$ is filled with $s$, we see that $s$ can not appear in the $(i+1)$-st row of $D(t\alpha)$.
Thus we have
\begin{align}\label{beta1}
\beta_1+\beta_2+\cdots+\beta_s\le (\alpha_1+\alpha_2+\cdots+\alpha_i)t\le (r_2+r_4+\cdots+r_{2k})t.
\end{align}
Moreover, since $s>r_1+r_2+\cdots+r_{2l-1}=d_{l-1}+r_{2l-1}\ge d_{l-1}$ and $\beta_s>0$, we derive that
\begin{align}\label{beta3}
\beta_1+\cdots+\beta_s>\beta_1+\cdots+\beta_{d_{l-1}}= (r_2+r_4+\cdots+r_{2l-2})t+c_1+\cdots+c_{l-1},
\end{align}
Combing \eqref{beta1} and \eqref{beta3}, we get
\[
(r_2+r_4+\cdots+r_{2l-2})t+c_1+\cdots+c_{l-1}
<(r_2+r_4+\cdots+r_{2k})t.
\]
Since $c_1+\cdots+c_{l-1}\ge0$, we must have $2l-2< 2k$, that is, $l<k+1$. This is a contradiction.
\qed

\section{Notched rectangle matroids}

In  this section, we consider Ehrhart polynomials of notched rectangle matroids.
We first explore some further properties of $F(a,b,c,t)$ as defined in \eqref{Fabct}, and then prove Corollary  \ref{coro2} and   Theorem \ref{spcia}.

By \eqref{eqj} and \eqref{origin} in the proof of Theorem \ref{main}, $F(a,b,c,t)$ is the number of integer solutions of the equation
\begin{align}\label{fc}
\left\{
  \begin{array}{ll}
    x_1+x_2+\cdots+x_{a+b}=bt+c, \\
    0\le x_i\le t,\ \text{for $1\le i\le a+b$}.
  \end{array}
\right.
\end{align}
We proceed to develop some further properties of the polynomial $F(a,b,c,t)$.

\begin{lem}\label{abba}
We have
\begin{align}
F(a,b,c,t)&=F(b,a,-c,t)\label{dh}\\
F(a,b,c,t)&=F(a+1,b-1,c+t,t)\label{shiftab}\\
F(a+1,b,0,t)&=\sum_{i=0}^tF(a,b,-i,t).\label{sfba}
\end{align}
\end{lem}

\pf Let $y_i=t-x_i$ in the equation \eqref{fc}. Then $y_1+y_2+\cdots+y_{a+b}=(a+b)t-(bt+c)=at-c$, where $0\le y_i\le t$. It is easy to see that the number of integer solutions of this equation is $F(b,a,-c,t)$.  Thus \eqref{dh} holds.

Since both  $F(a,b,c,t)$ and $F(a+1,b-1,t+c,t)$ are the number of solutions of $x_1+\cdots+x_{a+b}=bt+c=(b-1)t+t+c$, where $0\le x_i\le t$, we obtain \eqref{shiftab}.

Similarly, since $F(a+1,b,0,t)$ is the number of solutions of the equation $x_1+\cdots+x_{a+b}=bt-x_{a+b+1}$, where $0\le x_i\le t$ for $i=1,\ldots,a+b+1$, we see that \eqref{sfba} follows.  \qed

Let $a,b,c,d$ be nonnegative integers, and $r(S)=(a,b,c,d)$. Then  $u=(\min\{a,d\},0), v=(0,\min\{a,d\})$ and by \eqref{dh},
\begin{align}\label{generalf}
i((a,b,c,d),t)=i((d,c,b,a),t)=
\sum_{j=0}^{t\cdot\min\{a,d\}}F(a,b,j,t)F(c,d,-j,t).
\end{align}

\begin{theo}\label{m2sp}
For any nonnegative integers $a,b,c,d$, we have
\begin{align}\label{abcd}
i((a,b,c,d),t)+i((b,a,d,c),t)=i(U_{b+d,a+b+c+d},t)+i(U_{b,a+b},t)i(U_{d,c+d},t).
\end{align}
\end{theo}

\pf Since $F(a+c,b+d,0,t)$ is the number of solutions of $x_1+\cdots+x_{a+b+c+d}=(b+d)t$, where $0\le x_i\le t$, which is equal to the sum of number of solutions of
\begin{align}\label{dg}
\left\{
  \begin{array}{ll}
    x_1+\cdots+x_{a+b}=bt+j,\\
    x_{a+b+1}+\cdots+x_{a+b+c+d}=dt-j, \\
    0\le x_i\le t,\ \text{for $1\le i\le a+b+c+d$},
  \end{array}
\right.
\end{align}
for all possible integers $j$.
It is clear that if $j<-bt$, then the first equation of \eqref{dg} has no solutions. If $j<-ct$, then the second equation has no solution. Thus $j\ge-t\cdot\min\{b,c\}$. Similarly, one can check that $j\le t\cdot\min\{a,d\}$.
Therefore,
\begin{align*}
&F(a+c,b+d,0,t)\\&=
\sum_{j=-t\cdot\min\{b,c\}}^{t\cdot\min\{a,d\}}F(a,b,j,t)F(c,d,-j,t)\\[5pt]
&=\sum_{j=0}^{t\cdot\min\{a,d\}}F(a,b,j,t)F(c,d,-j,t) +\sum_{j=-t\cdot\min\{b,c\}}^{0}F(a,b,j,t)F(c,d,-j,t)\\
&\qquad -F(a,b,0,t)F(c,d,0,t)\\[5pt]
&=\sum_{j=0}^{t\cdot\min\{a,d\}}F(a,b,j,t)F(c,d,-j,t) +\sum_{j=0}^{t\cdot\min\{b,c\}}F(b,a,j,t)F(d,c,-j,t)\\
&\qquad -F(a,b,0,t)F(c,d,0,t)\\[5pt]
&=i((a,b,c,d),t)+i((b,a,d,c),t)-i(U_{b,a+b},t)i(U_{d,c+d},t),
\end{align*}
where the last step holds by \eqref{generalf}.  \qed

\noindent
{\bf Remark}. As pointed out by the referee, the statement of Theorem \ref{m2sp} is actually a fact of matroid subdivisions. More precisely, we split the uniform matroid polytope
\[
\P(U_{b+d,a+c+b+d})=\left \{ x\in\mathbb{R}_{\ge0} ^{a+b+c+d}:\sum_{i=1}^{a+b+c+d}x_i=b+d \right \}
\]
 into two pieces:
\[
\P_1=\left \{ x\in\mathbb{R}_{\ge0} ^{a+b+c+d}:\sum_{i=1}^{c+d}x_{a+b+i}\le d, \sum_{i=1}^{a+b+c+d}x_{i}=b+d \right \}
\]
and
\begin{align*}
\P_2=\left \{ x\in\mathbb{R}_{\ge0} ^{a+b+c+d}:\sum_{i=1}^{a+b}x_{i}\le b, \sum_{i=1}^{a+b+c+d}x_{i}=b+d\right\}.
%&=\left \{ x\in\mathbb{R} ^{a+b+c+d}:\sum_{i=1}^{c+d}x_{i}\ge d, \sum_{i=1}^{a+b+c+d}x_{i}=b+d\right\}.
\end{align*}
And their intersection is a common facet
\[
H=\left \{ x\in\mathbb{R}_{\ge0} ^{a+b+c+d}:\sum_{i=1}^{a+b}x_i=b, \sum_{i=1}^{c+d}x_{a+b+i}=d \right \}.
\]
One can check that
$\P_1=\P({\rm SM}_n(a, b, c, d))$ and $\P_2$ can be obtained by a rotation of the Schubert matroid polytope
\[
\P({\rm SM}_n(c,d,a,b))=\left \{ x\in\mathbb{R}_{\ge0} ^{a+b+c+d}:\sum_{i=1}^{a+b}x_{c+d+i}\le b, \sum_{i=1}^{a+b+c+d}x_{i}=b+d \right \}.
\]
Thus $i(\P_2,t)=i((c,d,a,b),t)=i((b,a,d,c),t)$. Moreover, it is easy to see that
 $i(H,t)=i(U_{b,a+b},t)i(U_{d,c+d},t)$.
Therefore,  we have
\begin{align*}
i(U_{b+d,a+b+c+d},t)=i((a,b,c,d),t)+i((b,a,d,c),t)
-i(U_{b,a+b},t)i(U_{d,c+d},t).
\end{align*}

\begin{coro}\label{xfa}
We have
\begin{align}\label{11ab}
i((1,1,a,b),t)+i((1,1,b-1,a+1),t)=(t+2)i(U_{b,a+b+1},t).
\end{align}
\end{coro}

\pf Since $F(1,1,i,t)=t+1-i$ for $0\le i\le t$, by \eqref{generalf}, we find
\begin{align*}
&i((1,1,a,b),t)+i((1,1,b-1,a+1),t)\\
&=\sum_{i=0}^t((t-i+1)F(a,b,-i,t)+(t-i+1)F(b-1,a+1,-i,t))\\
&=\sum_{i=0}^t((t-i+1)F(a,b,-i,t)+(t-i+1)F(a,b,i-t,t))\\
&=\sum_{i=0}^t((t-i+1)F(a,b,-i,t)+(i+1)F(a,b,-i,t))\\
&=\sum_{i=0}^t(t+2)F(a,b,-i,t)\\
&=(t+2)i(U_{b,a+b+1},t),
\end{align*}
where the second  step holds by \eqref{dh} and \eqref{shiftab}, and the last step follows from   \eqref{sfba}.
\qed

\noindent{\it Proof of Theorem \ref{spcia}.}
By \eqref{generalf}, we find
$i((a,b,a,b),t)=i((b,a,b,a),t)$.
The equations \eqref{ababs} and \eqref{aabb} follow directly from \eqref{abcd}. And \eqref{11aa1} is a special case of \eqref{11ab}. \qed

In the rest of this section, we give a proof of  Corollary \ref{coro2}.

Recall that the minimal matroid $T_{k,n}$ is the Schubert matroid ${\rm SM}_n(S)$ with
\[
S=\{2,3,\ldots,k,n\},
\]
where $n>k\ge2$.

\noindent{\it Proof of Corollary \ref{coro2}}. Since   $r=(1,k-1,n-k-1,1)$, by \eqref{generalf}, we have
\begin{align}
i(T_{k,n},t)=\sum_{j=0}^t F(1,k-1,j,t)F(n-k-1,1,-j,t).
\end{align}
By  \eqref{dh} and    \eqref{shiftab}, we have
\begin{align*}
F(1,k-1,j,t)=F(k-1,1,-j,t)=F(k,0,t-j,t)
=\binom{k+t-j-1}{t-j}
\end{align*}
and
\begin{align*}
F(n-k-1,1,-j,t)=F(n-k,0,t-j,t)={t-j+n-k-1\choose n-k-1}.
\end{align*}
Then
\begin{align}\label{ehrt}
i(T_{k,n},t)&=\sum_{j=0}^t{t-j+n-k-1\choose n-k-1}\binom{k+t-j-1}{t-j}\nonumber\\
&=\sum_{j=0}^t{j+n-k-1\choose n-k-1}\binom{k+j-1}{j}.
\end{align}
Thus we need to show that
\begin{align}\label{eq21}
\sum_{j=0}^t{j+n-k-1\choose n-k-1}\binom{k+j-1}{j}
=\frac{1}{{n-1\choose k-1}}{t+n-k\choose n-k}\sum_{j=0}^{k-1}{n-k+j-1\choose j}{t+j\choose j}.
\end{align}
Let $s=n-k$ in \eqref{eq21}, then we aim to show that
\begin{align}\label{feq}
\binom{s+k-1}{s} \sum_{j=0}^{t} \binom{j+s-1}{s-1} \binom{j+k-1}{k-1} =\binom{t+s}{s} \sum_{j=0}^{k-1} \binom{j+s-1}{s-1} \binom{t+j}{j}.
\end{align}
It is easy to see that the left hand side of \eqref{feq} is the coefficient of $y^{s-1}x^{k-1}$ in  \begin{align}\label{lfh}
&\binom{s+k-1}{s} \sum_{j=0}^{t} (1+y)^{s+j-1}(1+x)^{k+j-1}\nonumber\\&
 =\binom{k-1+s}{s}(1+y)^{s-1}(1+x)^{k-1}\cdot
 \frac{1-((1+y)(1+x))^{t+1}}{-x-y-xy}.
\end{align}
Similarly, the right hand side of  \eqref{feq} is the coefficient of $y^{s-1}x^{t}$ in
\begin{align}\label{rhs}
&\binom{t+s}{s} \sum_{j=0}^{k-1} (1+y)^{s+j-1}(1+x)^{t+j}\nonumber\\
 &=\binom{t+s}{s}(1+y)^{s-1}(1+x)^{t}\cdot
 \frac{1-((1+y)(1+x))^{k}}{-x-y-xy}.
 \end{align}
One can check that  the coefficient of $y^{s-1}x^{k-1}$ in \eqref{lfh} is equal to the coefficient of $y^{s-1}x^{t}$ in \eqref{rhs}. Thus \eqref{feq} follows. This completes the proof. \qed

To conclude this section, we remark that since $F(a,b,0,t)=i(U_{b,a+b},t)$ has  positive coefficients for any $a,b\ge1$, it is natural to ask whether $F(a,b,c,t)$ defined in \eqref{Fabct} has positive coefficients or not for any $c$.
The following conjecture  was verified for $a,b,c\le 10$.

\begin{conj}\label{cf}
$F(a,b,c,t)$  has positive coefficients for any $a,b\ge1$ if and only if $c=0,\pm1$.
\end{conj}

Since $F(1,1,c,t)=t+1-|c|$, we see that if $|c|>1$, then $F(1,1,c,t)$ has negative coefficients. Thus to prove Conjecture \ref{cf}, it is enough to show that if $c=\pm1$, then $F(a,b,c,t)$ is a positive polynomial in $t$ for any $a,b\ge1$.

\section{$(a,b)$-Catalan matroids}

In this section, we give a proof of Theorem \ref{coro3}.
Recall that the $(a,b)$-Catalan matroid $\mathbf{C}_n^{a,b}$ is the Schubert matroid ${\rm SM}_{(a+b)n}(S)$, where
\[r(S)=(\overbrace{a,b,a,b,\ldots,a,b}^{2n}).\]
Given a composition $\sigma$ of $n$,
$\ell(\sigma)$ denotes the number of parts of $\sigma$, $d(\sigma)$ denotes the cardinality of the equivalent class containing $\sigma$.
And $\Gamma_n$ is the set of pairwise non-equivalent compositions of $n$ with at least two parts and minimal parts larger than 1.

It can be seen  that \eqref{freef} is equivalent to
\begin{align}\label{pathm}
F(na,nb,0,t)&=n\cdot i(\mathbf{\overline{C}}_{n}^{a,b},t)+ i(U_{b,a+b},t)^n+\sum_{\sigma \in \Gamma _n} (-1)^{\ell(\sigma)-1}\frac{nd(\sigma )}{\ell(\sigma )}\cdot i(\overline{\mathbf{C}}_{\sigma}^{a,b},t),
\end{align}
where $\displaystyle i(\overline{\mathbf{C}}_{\sigma}^{a,b},t)=\prod_{j=1}^{\ell(\sigma )} i(\overline{\mathbf{C}}_{\sigma_j}^{a,b},t)$ and
\[
i(\overline{\mathbf{C}}_{\sigma_j}^{a,b},t)=i(\mathbf{C}_{\sigma_j}^{a,b},t)
-i(U_{b,a+b},t)\cdot i(\mathbf{C}_{\sigma_j-1}^{a,b},t).
\]
We shall prove \eqref{pathm} by  interpreting both sides in terms of weighted enumerations of certain lattice paths.

Let us begin with interpreting $i(\mathbf{C}_n^{a,b},t)$ and $F(na,nb,0,t)$ separately. By Theorem \ref{main}, since $r=(a,b,\ldots,a,b)$,  for $1\le j\le n$, we have
\begin{align}\label{uv}
u_j=\min\{a,(n-j)b\},\ v_j=\min\{b,(j-1)a\}.
\end{align}
Thus
\begin{align*}
i(\mathbf{C}_n^{a,b},t)&=\sum_{(c_1,\ldots,c_{n})}
\prod_{j=1}^{n}F(a,b,c_j,t),
\end{align*}
where $c_1+\cdots+c_{n}=0$,  $c_1+\cdots+c_j\ge0$   and $-tv_j\le c_j\le tu_j$,   for $1\le j\le n$.

On the other hand, since
$F(na,nb,0,t)$ is the number of  solutions of
\begin{align*}
\left\{
  \begin{array}{ll}
    x_1+x_2+\cdots+x_{(a+b)n}=bnt, \\
    0\le x_i\le t,\ \text{for $1\le i\le (a+b)n$}.
  \end{array}
\right.
\end{align*}
which is equivalent to the system of  equations
\begin{align}\label{fcsc}
\left\{
  \begin{array}{c}
 x_{1,1}+\cdots+x_{1,a+b}=bt+c'_1,   \\
  x_{2,1}+\cdots+x_{2,a+b}=bt+c'_2,   \\
    \vdots\\
 x_{n,1}+\cdots+x_{n,a+b}=bt+c'_n,
  \end{array}
\right.
\end{align}
for all possible integers $c'_1,\ldots,c'_n$,
where $0\le x_{i,j}\le t$ for $1\le i\le n$ and $1\le j\le a+b$, and $c'_1+\cdots+c'_{n}=0$.  It is easy to see that we can require  $-bt\le c'_j\le at$   for $1\le j\le n$. Thus
\begin{align*}
F(na,nb,0,t)=\sum_{(c'_1,\ldots,c'_{n})}
\prod_{j=1}^{n}F(a,b,c'_j,t),
\end{align*}
where $c'_1+\cdots+c'_{n}=0$,    and $-bt\le c'_j\le at$   for $1\le j\le n$.

Let
\begin{align}\label{pc}
\mathcal{C}_n^{a,b}=\left\{(c_1,\ldots,c_n)\,|\,\sum_{i=1}^nc_i=0,\  \sum_{i=1}^jc_i\ge0 \  \text{and}\  -tv_j\le c_j\le tu_j,\ \forall \ 1\le j\le n\right\}
\end{align}
and
\begin{align}\label{pfc}
\mathcal{F}_n^{a,b}=\left\{(c_1,\ldots,c_n)\,|\,\sum_{i=1}^nc_i=0,\    -bt\le c_j\le at,\ \forall \ 1\le j\le n\right\}.
\end{align}
Clearly, $\mathcal{C}_n^{a,b}\subseteq \mathcal{F}_n^{a,b}$.
We can view each sequence $\mathbf{c}=(c_1,\ldots,c_n)\in \mathcal{F}_n^{a,b}$ as a lattice path from $(0,0)$ to $(n,0)$ such that $c_j$ represents: an up step $(0,0)\rightarrow (1,c_j)$ if $c_j>0$, a down step $(0,0)\rightarrow (1,-|c_j|)$ if $c_j<0$, or a horizontal step $(0,0)\rightarrow(1,0)$ if $c_j=0$. Assign a weight to $\mathbf{c}$ as
\[{\rm wt}(\mathbf{c})=\prod_{j=1}^{n}F(a,b,c_j,t).\]
In particular, $F(a,b,0,t)^n$ is the weight of the path $(0,0,\ldots,0)$. Then
\[
i(\mathbf{C}_n^{a,b},t)
=\sum_{\c\in\mathcal{C}_n^{a,b}}{\rm wt}(\mathbf{c})\ \
\text{and}\ \ F(na,nb,0,t)
=\sum_{\c\in\mathcal{F}_n^{a,b}}{\rm wt}(\mathbf{c})
\]
can be viewed as weighted enumerations of lattice paths in $\mathcal{C}_n^{a,b}$ and $\mathcal{F}_n^{a,b}$, respectively.

Let
\begin{align}\label{barpc}
\mathcal{\overline{C}}_n^{a,b}=\{(c_1,\ldots,c_n)\in\mathcal{C}_n^{a,b}
\,|\,c_n\neq0\}.
\end{align}
Then we find
\begin{align}\label{barcehr}
i(\mathbf{\overline{C}}_n^{a,b},t)
=\sum_{\c\in\mathcal{\overline{C}}_n^{a,b}}{\rm wt}(\mathbf{c}).
\end{align}
Moreover, given a composition $\sigma=(\sigma_1,\ldots,\sigma_s)$ of $n$, denote
\begin{align}\label{barpc}
\mathcal{\overline{C}}_{\sigma}^{a,b}=
\{(\mathbf{c}_1,\ldots,\mathbf{c}_s)
\,|\,\mathbf{c}_j\in\mathcal{\overline{C}}_{\sigma_j}^{a,b},\ \text{for}\ 1\le j\le s\},
\end{align}
then
\[
i(\mathcal{\overline{C}}_{\sigma}^{a,b},t)=\prod_{j=1}^{\ell(\sigma )} i(\overline{\mathbf{C}}_{\sigma_j}^{a,b},t)
=\sum_{\c\in\mathcal{\overline{C}}_{\sigma}^{a,b}}{\rm wt}(\c).
\]
Therefore, \eqref{pathm} is equivalent to
\begin{align}\label{guod}
\sum_{\c\in\mathcal{F}_n^{a,b}}{\rm wt}(\mathbf{c})
=n\cdot\sum_{\c\in\mathcal{\overline{C}}_n^{a,b}}{\rm wt}(\mathbf{c})+{\rm wt}((0,\ldots,0))+\sum_{\sigma\in \Gamma_n}(-1)^{\ell(\sigma)-1}\frac{nd(\sigma )}{\ell(\sigma )}\sum_{\c\in\mathcal{\overline{C}}_{\sigma}^{a,b}}{\rm wt}(\c).
\end{align}

To prove \eqref{guod}, we shall enumerate the number of appearances of each $\c\in\mathcal{F}_n^{a,b}$ in the right hand side of \eqref{guod} by inclusion-exclusion. To this end, define two shifting operators $R$ and $L$ on each $\c=(c_1,\ldots,c_n)\in \mathcal{F}_n^{a,b}$ as follows
\begin{align*}
R(\c)&=(c_2,\ldots,c_n,c_1),\\
L(\c)&=(c_n,c_1,\ldots,c_{n-1}),
\end{align*}
and let $R^m(\c)=(c_{m+1},\ldots,c_n,c_1,\ldots,c_{m})$ denote the effect of applying $R$ to $\c$ $m$ times, where  $m\ge1$. Similar $L^m(\c)$ denote applying $L$ to $\c$ $m$ times.  It is clear that $R^n(\c)=L^n(c)=\c$.
Denote $\mathcal{L}^m_{\mathbf{c}}$ by the list of paths obtained from $\c$ by applying $L$ to $\c$ $m-1$ times, i.e.,
\[
\mathcal{L}^m_{\mathbf{c}}=(\c,L(\c),L^2(\c),\ldots,L^{m-1}(\c)).
\]
It is quite possible that $\mathcal{L}^m_{\mathbf{c}}$ contains repeated paths.
Similarly, let
\[
\mathcal{R}^m_{\mathbf{c}}
=(\c,R(\c),R^2(\c),\ldots,R^{m-1}(\c)).
\]
Obviously, each path in $\mathcal{L}^m_{\mathbf{c}}$ or $\mathcal{R}^m_{\mathbf{c}}$ has the same weight.

Denote $n\ast\mathcal{\overline{C}}_n^{a,b}$ by the list of paths that each  $\c$ of $\mathcal{\overline{C}}_n^{a,b}$ are replaced by the $n$ paths in $\mathcal{R}^n_{\mathbf{c}}$.
Similarly, denote $\frac{nd(\sigma )}{\ell(\sigma )}\ast\mathcal{\overline{C}}_{\sigma}^{a,b}$ by the list of paths that each  $\c$ of $\mathcal{\overline{C}}_{\sigma}^{a,b}$ are replaced by the   $\frac{nd(\sigma )}{\ell(\sigma )}$ paths in $\mathcal{R}^{\frac{nd(\sigma )}{\ell(\sigma )}}_{\mathbf{c}}$
for any $\sigma\in \Gamma_n$. If there is a minus sign in front of some $\frac{nd(\sigma )}{\ell(\sigma )}\ast\mathcal{\overline{C}}_{\sigma}^{a,b}$, then we delete the number of appearances of each $\c\in\frac{nd(\sigma )}{\ell(\sigma )}\ast\mathcal{\overline{C}}_{\sigma}^{a,b}$ in the enumeration. Assume that $\Gamma_n$ has $\gamma_n$ different compositions, i.e., $\Gamma_n=\{\sigma^1,\sigma^2,\ldots,\sigma^{\gamma_n}\}$.  Let
\begin{align}\label{jihe}
\left(n\ast\mathcal{\overline{C}}_n^{a,b},
(-1)^{\ell(\sigma^1)-1}\frac{nd(\sigma^1)}{\ell(\sigma^1 )}\ast\mathcal{\overline{C}}_{\sigma^{1}}^{a,b},
\ldots,(-1)^{\ell(\sigma^{\gamma_n})-1}
\frac{nd(\sigma^{\gamma_n})}{\ell(\sigma^{\gamma_n} )}\ast\mathcal{\overline{C}}_{\sigma^{\gamma_n}}^{a,b}
\right).
\end{align}
To prove \eqref{guod},
we aim to show  that, after cancellations,   each  $\c\in\mathcal{F}_n^{a,b}, \c\neq(0,\ldots,0)$  appears exactly once in \eqref{jihe}.

For example, let $n=6,a=2,b=3, t=1$ and  $\c=(1,-1,2,-2,1,-1)\in\mathcal{F}_6^{2,3}$. Then $\Gamma_6=\{(4,2),(3,3),(2,2,2)\}$. We aim to enumerate the number of appearances of $\c$ in
\begin{align}\label{lizz}
\left( 6\ast\overline{\mathcal{C}}_{6}^{2,3}, -6\ast \overline{\mathcal{C}}_{(4,2)}^{2,3},-3\ast \overline{\mathcal{C}} _{(3,3)}^{2,3},  2\ast\overline{\mathcal{C}}_{(2,2,2)}^{2,3}  \right).
\end{align}
One can check that there are $3$ appearances of $\c$   in $6\ast\overline{\mathcal{C}}_{6}^{2,3}$. That is, for the $3$ paths $\c_1=\c,\c_2=(2,-2,1,-1,1,-1),\c_3=(1,-1,1,-1,2,-2)$ in $\overline{\mathcal{C}}_{6}^{2,3}$, $\c$ appears in each of $\mathcal{R}^6_{\c_1}, \mathcal{R}^6_{\c_2},\mathcal{R}^6_{\c_3}$ exactly once. Similarly, $\c$
appears in $6\ast\overline{\mathcal{C}}_{(4,2)}^{2,3}$ $3$ times with a minus sign. That is, for the $3$ paths $\c_1'=\c,\c_2'=\c_2,\c_3'=\c_3$ in $\overline{\mathcal{C}}_{(4,2)}^{2,3}$, $\c$ appears in each of $\mathcal{R}^6_{\c_1'}, \mathcal{R}^6_{\c_2'},\mathcal{R}^6_{\c_3'}$ exactly once. Moreover, $\c$
appears in $3\ast\overline{\mathcal{C}} _{(3,3)}^{2,3}$ $0$ times, and appears in $2\ast\overline{\mathcal{C}}_{(2,2,2)}^{2,3}$ exactly once. That is,  $\c\in \overline{\mathcal{C}}_{(2,2,2)}^{2,3}$ and $\c$ only appears in $\mathcal{R}^2_\c=(\c,(-1,2,-2,1,-1,1))$ once. Therefore, the total number of appearances of $\c$ in \eqref{lizz} is $1$.
%\end{exam}

Given a path $\mathbf{c}\in\mathcal{F}_n^{a,b}$, if  $\mathbf{c}$ does not go below the $x$-axis, then  we write $\mathbf{c}\ge0$ for simplicity, and say $\mathbf{c}$ is \textit{nonnegative}. Otherwise,  write $\c<0$  and say   $\c$ is \textit{negative}.
Obviously, if $\c\ge0$, then $c_1+\cdots+c_j\ge0$ for any $1\le j\le n$.

It is clear that $\mathcal{\overline{C}}_n^{a,b}\subseteq\mathcal{F}_n^{a,b}$ and  $\mathcal{\overline{C}}_{\sigma}^{a,b}\subseteq\mathcal{F}_n^{a,b}$ for any $\sigma\in \Gamma_n$. By the definitions of $\mathcal{C}_n^{a,b}$ and $\mathcal{F}_n^{a,b}$ in \eqref{pc} and \eqref{pfc}, it seems possible that there exists $\c\in \mathcal{F}_n^{a,b}$ and $\c\ge0$, but $\c\notin\mathcal{C}_n^{a,b}$.
However, we show that this situation cannot happen.

For $\c=(c_1,\ldots,c_n)\in\mathcal{F}_n^{a,b}$, let
\[
\mathcal{E}_{\mathbf{c}}=\{(c_{j},\ldots,c_n,c_1,\ldots,c_{j-1})\,|\,1\le j\le n\}
\]
denote the set of paths that can be obtained from  $\c$ by cyclic shifting.

\begin{lem}\label{eqclass}
Each $\mathcal{E}_{\mathbf{c}}$ contains a nonnegative path.
\end{lem}

\pf
Given   $\mathbf{c}=(c_1,\ldots,c_n)$, define $d_j=c_1+\cdots+c_j$ for $1\le j\le n$. If $d_1,\ldots,d_n\ge0$, then it is clear that $\mathbf{c}$ is nonnegative. If   $d_j<0$ for some $j$, then let $k$ be the smallest index such that $d_k=\min\{d_1,\ldots,d_n\}$. Let $\mathbf{c}'=(c_{k+1},\ldots,c_n,c_1,\ldots,c_{k})\in \mathcal{E}_{\mathbf{c}}$. It is easy to see that $\mathbf{c}'$ is nonnegative.
\qed

\begin{lem}\label{fjf}
We have
\begin{align}\label{f2e}
\mathcal{F}_n^{a,b}=\bigcup_{\mathbf{c}\in
\mathcal{C}_n^{a,b}}
\mathcal{E}_{\mathbf{c}}.
\end{align}
\end{lem}
\pf
By   Lemma \ref{eqclass},
\[
\mathcal{F}_n^{a,b}=\bigcup_{\mathbf{c}\in
\mathcal{F}_n^{a,b},\, \mathbf{c}\ge0}
\mathcal{E}_{\mathbf{c}}.
\]
Clearly, $\mathcal{C}_n^{a,b}\subseteq \{\mathbf{c}\in\mathcal{F}_n^{a,b}\,|\, \mathbf{c}\ge0\}$. We aim to show that
\[
\{\mathbf{c}\in\mathcal{F}_n^{a,b}\,|\, \mathbf{c}\ge0\}\subseteq \mathcal{C}_n^{a,b}.
\]
Given $\mathbf{c}=(c_1,\ldots,c_n)\in\mathcal{F}_n^{a,b}$ such that $\mathbf{c}\ge0$, we need to show that $-tv_j\le c_j\le tu_j$ for any $1\le j\le n$. Since $\mathbf{c}\ge0$, we have $c_1+\cdots+c_{j-1}\ge0$. Adding the first $j$ equations in \eqref{fcsc} together, we obtain
\begin{align}\label{fcza}
jbt+c_1+\cdots+c_{j-1}+c_j=\sum_{i=1}^j\sum_{i'=1}^{a+b}x_{i,i'}\le bnt,
\end{align}
thus $c_j\le (n-j)bt$. Combing the fact $c_j\le at$, we arrive at
\[
c_j\le\min\{(n-j)bt,at\}=tu_j.
\]
On the other hand, by \eqref{fcza}, we obtain
\[
(j-1)bt+c_1+\cdots+c_{j-1}
=\sum_{i=1}^{j-1}\sum_{i'=1}^{a+b}x_{i,i'}\le (j-1)(a+b)t,
\]
thus
\[
(n-j+1)bt+c_j+\cdots+c_{n}
=\sum_{i=j}^{n}\sum_{i'=1}^{a+b}x_{i,i'}\ge bnt-(j-1)(a+b)t,
\]
so we see that
\begin{align}\label{cj}
c_j\ge-(j-1)at-(c_{j+1}+\cdots+c_n)\ge-(j-1)at,
\end{align}
where $c_{j+1}+\cdots+c_n\le0$ since $c_1+\cdots+c_n=0$ and $c_1+\cdots+c_j\ge0$. Combining the fact $c_j\ge-bt$ and \eqref{cj}, we obtain
\[
c_j\ge \min\{-bt,-(j-1)at\}=-tv_j,
\]
as desired. \qed

By Lemma \ref{eqclass} and Lemma \ref{fjf}, we can divide the discussions into three cases, that is, $\c\ge0, c_n\neq0$, or $\c\ge0, c_n=0$, or $\c<0$.

\begin{prop}\label{prop53}
Let $\c=(c_1,\ldots,c_n)\ge0$ such that  $c_n\neq0$. Then $\c$ appears exactly once in  \eqref{jihe}.
\end{prop}

To give a proof of Proposition \ref{prop53},  we need to enumerate how many copies of $\c$ appearing in each $\frac{nd(\sigma)}{\ell(\sigma)}\ast\mathcal{\overline{C}}_{\sigma}^{a,b}$. To this end, we first
 give   a combinatorial interpretation of the coefficient $\frac{nd(\sigma)}{\ell(\sigma)}$.

For a composition $\sigma=(\sigma_1,\ldots,\sigma_s)$ of $n$, denote $p(\sigma)$ by the least period of $\sigma$, that is, $p(\sigma)$ is the smallest integer such that $\sigma_i=\sigma_{i+p(\sigma)}$ for all $i$. If $\sigma$ has no period, then we define $p(\sigma)=\ell(\sigma)$. It is easy to see that $p(\sigma)=d(\sigma)$ is the cardinality of the equivalent class containing $\sigma$.   Let
\begin{align}
T(\sigma)=\sigma_1+\cdots+\sigma_{p(\sigma)}
\end{align}
be the sum of elements in a least period of $\sigma$.
Since $p(\sigma)=d(\sigma)$ and
$\frac{n}{T(\sigma)}=\frac{\ell(\sigma)}{p(\sigma)}$, we have
\begin{align}
T(\sigma)=\frac{nd(\sigma)}{\ell(\sigma)}.
\end{align}
For example, if $\sigma=(4,3,4,3,4,3,4,3)$, then $n=28, p(\sigma)=d(\sigma)=2, \ell(\sigma)=8,$ and $T(\sigma)=4+3=7$. If $\sigma=(3,3,2,2)$, then $n=10,p(\sigma)=d(\sigma)=\ell(\sigma)=4$ and $T(\sigma)=3+3+2+2=10$.

For $\c=(c_1,\ldots,c_n)\ge0$ and $c_n\neq0$,  we associate a unique composition
\begin{align}\label{picdf}
\pi(\c)=(\pi_1,\ldots,\pi_{\ell})
\end{align}
to $\c$ as follows. Assume that $\ell_1$ is the smallest index such that $c_1+\cdots+c_{\ell_1}=0$ and $c_{\ell_1}<0$. Then let $\pi_1=\ell_1$. Let $\ell_2$ be the smallest index such that $c_{\ell_1+1}+\cdots+c_{\ell_2}=0$ and $c_{\ell_2}<0$. Then let $\pi_2=\ell_2-\ell_1$. Continue this process, we can obtain  $\pi(\c)$ eventually. Since $\c\ge0$ and $c_n\neq0$, it is easy to see that $\pi_i\ge2$ for each $1\le i\le \ell$ and $\c\in\mathcal{\overline{C}}_{\pi(\c)}^{a,b}$. Moreover,
\begin{align}
\{\pi(\c)\,|\,\c\ge0,c_n\neq0\}=\Gamma_n\cup\{(n)\}.
\end{align}
For example, let $\c=(0,0,3,-1,0,-2,0,1,-1,0,1,1,-2,1,-1)$. Then $\pi(\c)=(6,3,4,2)$, see Figure \ref{pic} for an illustration of the path $\c$.

\begin{figure}[h]
\begin{center}
\begin{tikzpicture}[scale=0.4]

\draw [->][thick](-10mm,0mm)--(250mm,0mm);
\draw [->][thick](0mm,-5mm)--(0mm,45mm);
\draw[thick](0mm,0mm)--(15mm,0mm)--(30mm,0mm)
--(45mm,30mm)--(60mm,20mm)--(75mm,20mm)--(90mm,0mm)
--(105mm,0mm)--(120mm,10mm)--(135mm,0mm)--(150mm,0mm)
--(165mm,10mm)--(180mm,20mm)--(195mm,0mm)
--(210mm,10mm)--(225mm,0mm);

\node at (15mm,0mm) {$\cdot$};
\node at (165mm,10mm) {$\cdot$};

\node at (7mm,5mm) {\footnotesize $0$};
\node at (22mm,5mm) {\footnotesize $0$};
\node at (33mm,16mm) {\footnotesize $3$};
\node at (54mm,30mm) {\footnotesize $-1$};
\node at (67mm,25mm) {\footnotesize $0$};
\node at (86mm,15mm) {\footnotesize $-2$};
\node at (97mm,5mm) {\footnotesize $0$};
\node at (110mm,10mm) {\footnotesize $1$};
\node at (128mm,10mm) {\footnotesize $-1$};
\node at (142mm,5mm) {\footnotesize $0$};
\node at (157mm,11mm) {\footnotesize $1$};
\node at (171mm,20mm) {\footnotesize $1$};
\node at (190mm,15mm) {\footnotesize $-2$};
\node at (201mm,9mm) {\footnotesize $1$};
\node at (221mm,9mm) {\footnotesize $-1$};

\end{tikzpicture}

\end{center}
\vspace{-5mm}
\caption{The path $\c$.}\label{pic}
\end{figure}

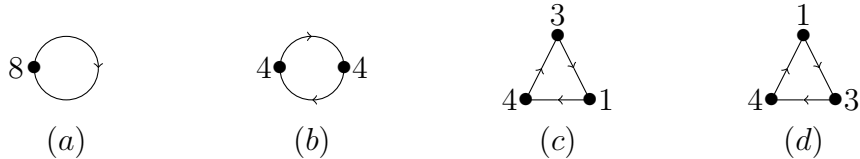
\begin{figure}[h]
\begin{center}
\begin{tikzpicture}[scale=0.4]
\draw [->] (0,0) arc (180:0:30pt);
\node at (0,0) {$\bullet$};\node at (-.6,0) {8};
\draw [-] (60pt,0) arc (0:-180:30pt);

\draw [->] (230pt,0pt) arc (180:90:30pt);
\draw [->] (260pt,30pt) arc (90:0:30pt);
\node at (230pt,0pt) {$\bullet$};
\node at (215pt,0pt) {4};\node at (305pt,0pt) {4};
\draw [->] (290pt,0pt) arc (0:-90:30pt);
\draw [-] (260pt,-30pt) arc (-90:-180:30pt);
\node at (290pt,0pt) {$\bullet$};

\node at (460pt,-30pt) {$\bullet$};
\node at (445pt,-30pt) {$4$};
\node at (520pt,-30pt) {$\bullet$};
\node at (535pt,-30pt) {$1$};
\node at (490pt,30pt) {$\bullet$};
\node at (490pt,50pt) {$3$};

\draw [->] (460pt,-30pt)--(475pt,0pt);
\draw (475pt,0pt)--(490pt,30pt);

\draw [->] (490pt,30pt)--(505pt,0pt);
\draw (505pt,0pt)--(520pt,-30pt);

\draw [->]  (520pt,-30pt)--(490pt,-30pt);
\draw (490pt,-30pt)--(460pt,-30pt);

%%%%%%%%%%%%%

\node at (690pt,-30pt) {$\bullet$};
\node at (675pt,-30pt) {$4$};
\node at (750pt,-30pt) {$\bullet$};
\node at (765pt,-30pt) {$3$};
\node at (720pt,30pt) {$\bullet$};
\node at (720pt,50pt) {$1$};

\draw [->] (690pt,-30pt)--(705pt,0pt);
\draw (705pt,0pt)--(720pt,30pt);

\draw [->] (720pt,30pt)--(735pt,0pt);
\draw (735pt,0pt)--(750pt,-30pt);

\draw [->]  (750pt,-30pt)--(720pt,-30pt);
\draw (720pt,-30pt)--(690pt,-30pt);

\node at (30pt,-70pt) {$(a)$};
\node at (260pt,-70pt) {$(b)$};
\node at (490pt,-70pt) {$(c)$};
\node at (720pt,-70pt) {$(d)$};

\end{tikzpicture}

\end{center}
\vspace{-5mm}
\caption{Circle representations of $(8),(4,4),(4,3,1),(4,1,3)$.}\label{cycg}
\end{figure}

Given a composition $\sigma=(\sigma_1,\ldots,\sigma_s)$, arrange $\sigma_1,\ldots,\sigma_s$ on a directed circle, such that there is a directed edge from $\sigma_i$ to $\sigma_{i+1}$ for $1\le i\le s-1$, and a directed edge from $\sigma_s$ to $\sigma_1$.  If $\sigma=(n)$ has only one part, then there is a directed loop on the node $n$. We call such a configuration the \textit{circle representation} of $\sigma$, denoted as $G(\sigma)$.  We view all the edges in $G(\sigma)$  different, even if they have the same nodes and directed edges.  For example, Figure \ref{cycg} displays the circle representations of $(8),(4,4),(4,3,1),(4,1,3)$, respectively. There are two different edges in Figure \ref{cycg}(b).

By \textit{contracting} a directed edge, say $\sigma_i\rightarrow\sigma_{i+1}$, of  $G(\sigma)$, we mean delete this edge and form a new node labeled by $\sigma_i+\sigma_{i+1}$, and keep all the other edges unchanged. Since all the edges in $G(\sigma)$ are viewed  different, it is quite possible that different ways of contracting the edges lead to the same circle representation.
For example, Figure \ref{cycg}(b) can be viewed as obtained by contracting the directed edge $3\rightarrow 1$ in Figure \ref{cycg}(c), or contracting the directed edge $1\rightarrow 3$ in Figure \ref{cycg}(d). Moreover, we can obtain Figure \ref{cycg}(a) twice by contracting the two different edges of Figure \ref{cycg}(b).

If $\sigma$ and $\sigma'$  have the same circle representation, then $\sigma\sim\sigma'$ are equivalent. It is easy to see that $T(\sigma)\ast
\mathcal{\overline{C}}_{\sigma}^{a,b}$ and  $T(\sigma')\ast
\mathcal{\overline{C}}_{\sigma'}^{a,b}$ contain the same number of appearances of each $\c\in\mathcal{F}_n^{a,b}$. Given a composition $\tau$, after contracting some edges of $G(\tau)$, we obtain a new circle representation, which is $G(\sigma)$ for some composition $\sigma$. To read off a specific $\sigma$, we can choose any node in $G(\sigma)$ as the first element $\sigma_1$, and then read off $\sigma_2,\sigma_3$, etc. of $\sigma$ from $G(\sigma)$ clock-wisely.

\begin{prop}\label{lemma55}
Let $\c\ge0$ and $c_n\neq0$. Then the number of  appearances of $\c$  in $T(\sigma)\ast\mathcal{\overline{C}}_{\sigma}^{a,b}$ is equal to the number of ways that $G(\sigma)$ can be obtained   by contracting edges in $G(\pi(\c))$.
\end{prop}

\pf
Let $\c=(c_1,\ldots,c_n)$ and $\pi(\c)=(\pi_1,\ldots,\pi_{\ell})$. We first show that if there is a way of contracting   edges of $G(\pi(\c))$ to obtain $G(\sigma)$, then $\c$ appears at least once in $T(\sigma)\ast\mathcal{\overline{C}}_{\sigma}^{a,b}$. Then we show that if $\c$ appears once in $T(\sigma)\ast\mathcal{\overline{C}}_{\sigma}^{a,b}$, then there is a way of contracting edges of $G(\pi(\c))$ to obtain $G(\sigma)$.

A contracting of edges   of $G(\pi(\c))$ is equivalent to adding consecutive elements of $\pi(\c)$ together, where we arrange $\pi(\c)$ on a circle, thus $\pi_1$ and  $\pi_{\ell}$ can be added together.
After contracting edges of  $G(\pi(\c))$, we obtain $G(\sigma)$.  Since $\sigma$ may have a period, to read off $\sigma$,  we need to locate a position of $\sigma_1$, and then read off $\sigma_2,\sigma_3,$ etc. from $G(\sigma)$ clock-wisely.  There are two cases, depending on whether $\pi_1$ and $\pi_{\ell}$ are added together or not.

Case 1. There exist $1\le j<i\le {\ell}$ such that $\pi_i,\ldots,\pi_{\ell},\pi_1,\ldots,\pi_j$ are added together. Let
$\sigma_1=\pi_i+\cdots+\pi_{\ell}+\pi_1+\cdots+\pi_j$.

Case 2. There exists $1\le i\le \ell$ such that $\pi_i,\ldots,\pi_\ell$ are added together. Let
$\sigma_1=\pi_i+\cdots+\pi_\ell$.

For both cases,  let
\[\c'=L^{\pi_i+\cdots+\pi_{\ell}}(\c).\]
One can check that $\c'\in\mathcal{\overline{C}}_{\sigma}^{a,b}$.
Since $\pi_i+\cdots+\pi_{\ell}\le T(\sigma)$ and $\c=R^{\pi_i+\cdots+\pi_{\ell}}(\c')$, we find that $\c$ will appear in $\mathcal{R}^{T(\sigma)}_{\c'}=(\c',R(\c'),\ldots,
R^{T(\sigma)-1}(\c'))$ at least once.

In the following, we show that if $\c$ appears in $T(\sigma)\ast\mathcal{\overline{C}}_{\sigma}^{a,b}$ once, then there is a way of contracting edges of $G(\pi(\c))$ to obtain $G(\sigma)$.

Suppose that there exists $\c'\in\mathcal{\overline{C}}_{\sigma}^{a,b}$ such that $\mathcal{R}^{T(\sigma)}_{\c'}=(\c',R(\c'),\ldots,
R^{T(\sigma)-1}(\c'))$ contains $k_0$ copies of $\c$. We aim to construct $k_0$ different ways of contracting edges of $G(\pi(\c))$ to obtain $G(\sigma)$.
Let $0\le i_1<i_2<\cdots<i_{k_0}\le T(\sigma)-1$ such that
\[R^{i_1}(\c')=R^{i_2}(\c')=\cdots=R^{i_{k_0}}(\c')=\c.\]
Then $\c'=L^{i_1}(\c)=\cdots=L^{i_{k_0}}(\c)$ and there exist $i_1'<i_2'<\cdots<i_{k_0}'$ such that
\[R^{i_1'}(\pi(\c'))=R^{i_2'}(\pi(\c'))=\cdots
=R^{i_{k_0}'}(\pi(\c'))=\pi(\c).\]
Since $\c'\in\mathcal{\overline{C}}_{\sigma}^{a,b}$, we can
add consecutive elements of $\pi(\c')=(\pi_1',\ldots,\pi_z')$ to obtain $\sigma=(\sigma_1,\ldots,\sigma_s)$. If we require that $\pi_1'$ and $\pi_z'$ can not be added together, then there are integers $j_1<j_2<\cdots<j_s=z$ such that
\[\sigma_1=\pi_1'+\cdots+\pi_{j_1}', \sigma_2=\pi_{j_1+1}'+\cdots+\pi_{j_2}',\ldots, \sigma_s=\pi_{j_{s-1}+1}'+\cdots+\pi_{j_s}'.\]
For the $k_0$ appearances of $\c$ in $\mathcal{R}^{T(\sigma)}_{\c'}$, we can construct $k_0$ ways of contracting edges of $G(\pi(\c))$ as following. For each $j\in\{i_1',\ldots,i_{k_0}'\}$, we can add the elements in $L^j(\pi(\c))$ with the same positions of elements in $\pi(\c')$. More precisely, for $j\in\{i_1',\ldots,i_{k_0}'\}$, let
\[\sigma_1=\pi_{1-j}+\cdots+\pi_{j_1-j},
\sigma_2= \pi_{j_1+1-j}+\cdots+\pi_{j_2-j},\ldots, \sigma_s=\pi_{j_{s-1}+1-j}+\cdots+\pi_{j_s-j},\]
where   the indices are taken  modulo $j_s=z$. It is easy to see that these $k_0$ constructions correspond to $k_0$ different ways of contracting edges of $G(\pi(\c))$.

For example, let $\pi(\c)=(2,3,2,2,3,2,2,3,2,2,3,2)$, $\sigma=(14,14)$ and $\pi(\c')=(2,2,3,$ $2,2,3,2,2,3,2,2,3)$. Then $i_1=2, i_2=9$, i.e., $ R^2(\c')=R^9(\c')=\c$ and $i_1'=1,i_2'=4$. Since
$j_1=6,j_2=12$, for $j=i_1'=1$,  let
$\sigma_1=\pi_{12}+\pi_1+\cdots+\pi_{5} $ and $\sigma_2=\pi_{6}+
\cdots+\pi_{11}$. For $j=i_2'=4$,  let
$\sigma_1=\pi_9+\cdots+\pi_{12}+\pi_1+\pi_{2} $ and $
\sigma_2=\pi_{3}+
\cdots+\pi_{8}$.

Moreover, if there is another $\c''\in\mathcal{\overline{C}}_{\sigma}^{a,b}$ such that $\mathcal{R}^{T(\sigma)}_{\c''}=(\c'',R(\c''),\ldots,
R^{T(\sigma)-1}(\c''))$ contains $m_0$ copies of $\c$.
Let \[R^{i_1'}(\c'')=R^{i_2'}(\c'')=\cdots=R^{i_{m_0}'}(\c'')=\c.\]
Then we must have $\{i_1,\ldots,i_{k_0}\} \cap\{i_1',\ldots,i_{m_0}'\}=\emptyset$.
In fact, if there exists $j\in\{i_1,\ldots,i_{k_0}\}\cap\{i_1',\ldots,i_{m_0}'\}$, then $R^j(\c')=R^j(\c'')=\c$, this implies $\c'=\c''$. By the same constructions for $\c'$, we can obtain $m_0$  ways of contracting edges of $G(\pi(\c))$, different from the above $k_0$ ways. This completes the proof. \qed

\noindent{\it Proof of Proposition \ref{prop53}}.
Assume that $\pi(\c)$ has $m=\ell(\pi(\c))$ parts, i.e., $G(\pi(\c))$ has $m$ edges. It is clear that the number of appearances of $\c$   in $(-1)^{m-1}T(\pi(\c))\ast\mathcal{\overline{C}}_{\pi(\c)}^{a,b}$ is $(-1)^{m-1}$, corresponding to contracting 0 edges in   $G(\pi(\c))$.
By Proposition \ref{lemma55}, if we contract   one edge of $G(\pi(\c))$, then $\c$ will be enumerated by $m$ times, with sign $(-1)^{m}$. Similarly, if we contract any two edges of $G(\pi(\c))$, $\c$ will be enumerated by ${m\choose 2}$ times with sign $(-1)^{m+1}$, etc. Therefore, the total number of appearances of $\c$ is
\[
(-1)^{m-1}{m\choose0}
+(-1)^{m}{m\choose1}
+(-1)^{m+1}{m\choose2}
+\cdots
+(-1)^{m+m-2}{m\choose m-1}=1,
\]
as required. \qed

For the running example in  \eqref{lizz}, given $\c=(1,-1,2,-2,1,-1)$, we have $\pi(\c)=(2,2,2)$. We need to contract  edges of $G((2,2,2))$ to obtain circle representations of   $(6),(4,2),(3,3),(2,2,2)$.
For $\sigma=(6)$, there are 3 ways of contracting 2 edges among all 3 edges in $G((2,2,2))$ to obtain $G((6))$, so    $\c$ is counted $\binom{3}{2}$  times in $6\ast\overline{C} ^{2,3}_{(6)}$. For $\sigma=(4,2)$, there are 3 ways to  contract one edge of $G((2,2,2))$ to obtain $G((4,2))$, so $\c$ is counted $\binom{3}{1}$  times in $6\ast\overline{C} ^{2,3}_{(4,2)}$ with a minus sign. Similarly, we can contract 0 edges in $G((2,2,2))$ to obtain $G((2,2,2))$, which means $\c$ is counted $\binom{3}{0}$  times in $2\ast\overline{C} ^{2,3}_{(2,2,2)}$. We can not obtain $G((3,3))$ by contracting any edges of $G((2,2,2))$.   Consequently, the total number of appearances of $\c$ in \eqref{lizz} is
\[\binom{3}{2}-\binom{3}{1}+\binom{3}{0}=1.\]

In the following, we consider the case $\c=(c_1,\ldots,c_n)<0$.

\begin{prop}\label{pp56}
If $\c<0$, then    $\c$  appears in \eqref{jihe} exactly once.
\end{prop}

\pf  For $1\le j\le n$, let $d_j=c_1+\cdots+c_j$. Let $k$ be the largest index such that $d_k=\min\{d_1,\ldots,d_n\}$.  Define
\[
\c'=L^{n-k}(\c).
\]
It is easy to see that $\c'\in\mathcal{\overline{C}}_n^{a,b}$ and $n-k$ is the smallest integer $m$ such that $L^{m}(\c)\ge0$. By Proposition \ref{prop53}, the total number of appearances of $\c'$ in \eqref{jihe} is 1.
 For example, let $\c=(2,-4,7,-7,4,-4,4,-2)$. Then $k=6$ and $\c'=(4,-2,2,-4,7,-7,4,-4)$. The paths $\c$ and $\c'$ are displayed in Figure  \ref{ccbar}.

\begin{figure}[h]
\begin{center}
\begin{tikzpicture}[scale=0.4]

\draw [->][thick](0mm,0mm)--(100mm,0mm);
\draw [->][thick](5mm,-5mm)--(5mm,30mm);
\draw[thick](5mm,0mm)--(15mm,10mm)
--(25mm,-10mm)--(35mm,25mm)--(45mm,-10mm)
--(55mm,10mm)--(65mm,-10mm);
\draw [ultra thick](65mm,-10mm)--(75mm,10mm)
--(85mm,0mm);

\draw [dashed] (65mm,0mm)--(65mm,-10mm);
\node at (65mm,5mm) {$k$};
\draw [dashed] (175mm,0mm)--(175mm,10mm);
\node at (175mm,-5mm) {$n$$-$$k$};
\node at (235mm,-5mm) {$n$};
\node at (85mm,-5mm) {$n$};
\node at (40mm,-20mm) {$\c$};
\node at (190mm,-20mm) {$\c'$};

\draw [->][thick](150mm,0mm)--(250mm,0mm);
\draw [->][thick](155mm,-5mm)--(155mm,30mm);
\draw[ultra thick](155mm,0mm)--(165mm,20mm)
--(175mm,10mm);
\draw[thick](175mm,10mm)--(185mm,20mm)--(195mm,0mm)
--(205mm,35mm)--(215mm,0mm)--(225mm,20mm)
--(235mm,0mm);

\end{tikzpicture}

\end{center}
\vspace{-5mm}
\caption{The paths $\c$ and $\c'$.}\label{ccbar}
\end{figure}

For any path $\c''\in\mathcal{\overline{C}}_{\sigma}^{a,b}$, let
\[
\mathcal{R}^{T(\sigma)}(\c'')=(\c'',R(\c''),R^2(\c''),
\ldots,R^{T(\sigma)-1}(\c'')).
\]
We aim to show that the number of appearances of $\c$ and $\c'$ in $\mathcal{R}^{T(\sigma)}(\c'')$ are the same. We first show that if there are two paths $\c$   in $\mathcal{R}^{T(\sigma)}(\c'')$, then there is a path $\c'$ between them. Suppose that there exist $i<j$ such that $R^i(\c'')=R^j(\c'')=\c$.  Then
\begin{align}\label{ccp}
\c'=L^{n-k}(\c)=L^{n-k}(R^j(\c''))=R^{j-(n-k)}(\c'').
\end{align}
Since $j-i>n-k$, we have $i<j-(n-k)<j$, which means that $\c'$ appears between $R^i(\c'')$ and $R^j(\c'')$ at least once.
Similarly, we can obtain that there is a $\c$ between any two  $\c'$ in $\mathcal{R}^{T(\sigma)}(\c'')$.

At this moment, we can only conclude that the number of appearances of $\c$ and $\c'$ in $\mathcal{R}^{T(\sigma)}(\c'')$ are equal or differ  by 1. Let $i_0$ be the smallest index such that $R^{i_0}(\c'')=\c$ and $j_0$ be the largest index such that $R^{j_0}(\c'')=\c$.  Since $\c''=L^{i_0}(\c)\ge0$, we find that $i_0\ge n-k$. Since $\c'=L^{n-k}(\c)$, we find that $\c'$ must appear in $(\c'',R(\c''),\ldots,R^{i_0}(\c''))$.
Similarly, we can show that there is a $\c$ appearing to the right of the right-most $\c'$. This completes the proof.
\qed

Finally, we consider the case $\c\ge0, c_n=0$.

\begin{prop}
 If $\c\ge0, c_n=0$, then    $\c$ appears in \eqref{jihe} exactly once.
\end{prop}

\pf
Let $q$ be the largest index such that $c_q\neq0$ and $c_{q+1}=\cdots=c_n=0$.  Define
\[
\c'=L^{n-q}(\c)=(c_{q+1},\ldots,c_n,c_1,\ldots,c_q).
\]
Similar to the proof   of Proposition \ref{pp56}, we can show that the number of appearances of $\c$ and $\c'$   are exactly the same in each $T(\sigma)\ast\mathcal{\overline{C}}_{\sigma}^{a,b}$. \qed

\section{Sparse paving Schubert matroids}

In this section, we study  sparse paving Schubert matroids. Recall that a matroid $M$ if sparse paving if and only if every subset of cardinality ${\rm rk}(M)$ is either a basis or a circuit-hyperplane.

Let $r=(k-1,1,1,n-k-1)$, where $n>k\ge2$, we obtain a Schubert matroid, denoted as ${\rm Sp}_{k,n}$. As will be shown, ${\rm Sp}_{k,n}$ is sparse paving, and a Schubert matroid is sparse paving if it is either a uniform matroid or a uniform matroid with one basis removed. Equivalently, a sparse paving Schubert matroid  is a notch matroid  with the upper bounding path $U=N^{n-k}E^{k}$, see Bon and de Mier \cite[Definition 8.1]{Bon2}.

\begin{prop}\label{sps}
A  Schubert matroid ${\rm SM}_n(S)$ of rank $n-k$ is sparse paving if and only if it is uniform or $r(S)=(k-1,1,1,n-k-1)$, i.e.,
\[
S=\{k+1,\ldots,n\}\ \ \text{or}\ \ \{k,k+2,\ldots,n\}.
\]
\end{prop}

\pf If  $S=\{k+1,\ldots,n\}$, then ${\rm SM}_n(S)$ is the uniform matroid $U_{n-k,n}$, which is sparse paving by definition. Suppose that  $S=\{k,k+2,\ldots,n\}$. Then we aim to show that every $(n-k)$-subset  of $[n]$ is either a basis or a circuit-hyperplane. In fact, it is easy to see that there is exactly one $(n-k)$-subset of $[n]$ which is not a basis of ${\rm SM}_n(S)$, i.e., $T=\{k+1,\ldots,n\}$. It is also straightforward to check that $T$ is both a circuit and a flat. Moreover, $T$ is also a hyperplane since ${\rm rk}_S(T)=n-k-1$. Thus ${\rm SM}_n(S)$ is sparse paving.

On the contrary, suppose that a Schubert matroid $M={\rm SM}_n(S)$ is sparse paving and it is not uniform. Then every subset of cardinality $n-k$ is either a basis or a circuit-hyperplane. Since $M$ is not uniform, it has a circuit-hyperplane. It is easy to see that $M$ is connected, by \cite[Theorem 8.3]{Bon2}, $M$ is a notch matroid. Thus $r(S)=(k-1,1,1,n-k-1)$. \qed

By the proof of Proposition \ref{sps}, if $S=\{k,k+2,\ldots,n\}$, then its corresponding sparse paving Schubert matroid ${\rm Sp}_{k,n}$ has exactly one circuit-hyperplane, i.e., $\{k+1,\ldots,n\}$. The following lemma follows from Ferroni \cite[Corollary 4.6]{Fer3}.

\begin{lem}\label{umk}
We have
\[
i({\rm Sp}_{k,n},t)=i({\rm Sp}_{n-k,n},t)=i(U_{k,n},t)-i(T_{k,n},t-1).
\]
\end{lem}

By Lemma \ref{umk}, we find that $i(U_{k,n},t)-i({\rm Sp}_{k,n},t)=i(T_{k,n},t-1)$.
By Ferroni \cite[Theorem 1.6]{Fer2}, $i(T_{k,n},t-1)$ has positive coefficients. Therefore,
\[
i({\rm Sp}_{k,n},t)\le i(U_{k,n},t).
\]
To finish the proof of Theorem \ref{spp},
we still need  to show that $i(T_{k,n},t)\le i({\rm Sp}_{k,n},t)$, that is,
\[
i(U_{k,n},t)-i(T_{k,n},t-1)-i(T_{k,n},t)
\]
has positive coefficients.
We aim to prove two slightly stronger statements, i.e.,
\begin{align}\label{u2k}
i(U_{k,n},t)\ge i(U_{2,n},t)
\end{align}
and
\begin{align}\label{utt}
i(U_{2,n},t)-i(T_{k,n},t-1)-i(T_{k,n},t)\ge0.
\end{align}

In order to prove \eqref{u2k}, we need the combinatorial interpretation of the coefficients of $i(U_{k,n},t)$.  For any $m\ge0$,
Ferroni \cite[Theorem 4.3]{Fer} showed that the coefficient of $t^m$ in $i(U_{k,n},t)$ is
\begin{align}\label{ukn}
[t^m]i(U_{k,n},t)=\frac{1}{(n-1)!}\sum_{j=0}^{k-1}W(j,n,m+1)A(m,k-j-1),
\end{align}
where $W(j,n,m+1)$ are the weighted Lah numbers and $A(m,k-j-1)$ are the Eulerian numbers. In particular,
$W(0,n,k)={n\brack k}$
is the unsigned Stirling number of the first kind. The following properties of ${n\brack k}$ are well  known,
\[
t(t+1)\cdots (t+n-1)=\sum_{k=0}^n{n\brack k}t^k,
\]
and
\begin{align}\label{stirling}
{n\brack k}=(n-1){n-1\brack k}+{n-1\brack k-1},
\end{align}
see, for example, Stanley \cite{sta}.

\begin{lem}\label{uk2}
For $3\le k\le\frac{n}{2}$, we have the coefficient-wise relation
\[i(U_{2,n},t)\le i(U_{k,n},t).\]
\end{lem}

\pf
By \eqref{ukn},  we have
\[[t^m]i(U_{2,n},t)=\frac{1}{(n-1)!}(W(0,n,m+1)A(m,1)
+W(1,n,m+1)A(m,0)).
\]
There are four cases to consider.

Case 1. $k<m$. By \eqref{ukn},
\[[t^m]i(U_{k,n},t)\ge\frac{1}{(n-1)!}(W(0,n,m+1)A(m,k-1)
+W(1,n,m+1)A(m,k-2)).\]
Since the Eulerian numbers are unimodal and $k-1\ge2$, we have   $A(m,k-1)\ge A(m,1)$ and $A(m,k-2)\ge A(m,0)$. So
$[t^m]i(U_{k,n},t)\ge[t^m]i(U_{2,n},t).$

Case 2. $k=m=3$.  When $n=4,5$, the lemma holds obviously.
When $n\ge6$, it is easy to see that
\[W(0,n,4)\le W(1,n,4)\ \text{and} \
W(1,n,4)\le W(2,n,4).\]
By \eqref{ukn},
\begin{align*}
[t^3]i(U_{3,n},t)
&=\frac{1}{(n-1)!}(W(0,n,4)
+4W(1,n,4)+W(2,n,4))\\
&\ge\frac{1}{(n-1)!}(
4W(1,n,4)+W(2,n,4))\\
&\ge\frac{1}{(n-1)!}(4W(0,n,4)
+W(1,n,4)))\\
&=[t^3]i(U_{2,n},t).
\end{align*}

Case 3. $k=m\ge4$.     By \eqref{ukn},
\begin{align*}
[t^m]i(U_{k,n},t)
&\ge\frac{1}{(n-1)!}(W(1,n,m+1)A(m,m-2)
+W(m-2,n,m+1)A(m,1)).
\end{align*}
 Since
$A(m,m-2)\ge A(m,0)=1,$   it is enough to show that
\begin{align}\label{ww}
W(m-2,n,m+1)\ge W(0,n,m+1).
\end{align}
We shall give a combinatorial proof of \eqref{ww}.

 Denote $\mathcal{W}(\ell, n, m)$ by the set of partitions of weight $\ell$ of $[n]$     into   $m$ linearly ordered blocks, and let $W(\ell, n, m)$ be the cardinality of $\mathcal{W}(\ell, n, m)$, see Ferroni \cite{Fer}. For  a partition $\pi$  with linearly ordered blocks, the weight of $\pi$ is
$w(\pi)=\sum_{B\in \pi}w(B)$,
where $w(B)$ is the number of elements in  $B$ that are smaller than the first element in $B$. In order to prove \eqref{ww}, we construct an injection from $\mathcal{W}(0, n, m+1)$ to $\mathcal{W}(m-2, n, m+1)$.

By definition,
$\mathcal{W}(0, n, m+1)$ is the set of partitions   of $[n]$     into   $m+1$  blocks, the elements of each block are arranged increasingly. Suppose that $\tau\in\mathcal{W}(0, n, m+1)$. We aim to construct a partition $\tau'\in\mathcal{W}(m-2, n, m+1)$ from $\tau$. Let $B_1,\ldots,B_j$ be the blocks of $\tau$ having more than one element and the smallest element of $B_i$ is smaller than the smallest element of $B_{i+1}$ for $1\le i\le j-1$. Apparently,
\[|B_1|+\cdots+|B_j|=n-(m+1-j).\]
For any linearly ordered  block $B=(b_1,\ldots,b_s)$ with $w(B)=0$, that is, $b_1$ is the smallest element of $B$. Assume that $b_m$ is the largest of element of $B$. Let $B'$ be obtained from $B$ by cyclically shifting $b_m$ to the first position, i.e., $B'=(b_m,\ldots,b_s,b_1,\ldots,b_{m-1})$. Then  $w(B')=|B|-1$. Therefore,
\begin{align*}
w(B_1')+\cdots+w(B_j')&=|B_1|-1+\cdots+|B_j|-1=n-(m+1).
\end{align*}

Since $k=m\le\frac{n}{2}$, we have $n-(m+1)\ge m-1$. Thus we can construct a partition $\tau'\in\mathcal{W}(m-2, n, m+1)$ with weight $m-2$ from  $\tau$ as follows. There exists some index $i\, (1\le i<j)$ such that
\[
w(B_1')+\cdots+w(B_i')\le m-2\ \text{and}\ w(B_1')+\cdots+w(B_{i+1}')>m-2.
\]
We can cyclically shift a suitable element of $B_{i+1}$ to  the first position to obtain a new block $B_{i+1}''$, such that \[w(B_1')+\cdots+w(B_{i}')+w(B_{i+1}'')=m-2.\]
Keep the other blocks of $\tau'$ the same with those of $\tau$. It is easy to see that this construction is an injection from  $\mathcal{W}(0, n, m+1)$ to $\mathcal{W}(m-2, n, m+1)$.   This completes the proof of \eqref{ww}.

Case 4. $k>m$. By \eqref{ukn},
\[[t^m]i(U_{k,n},t)
\ge\frac{1}{(n-1)!}(W(k-1,n,m+1)A(m,0)
+W(k-2,n,m+1)A(m,1)).\]
We aim to show that
\begin{align}
\varphi (n,k,m):=W(k-2,n,m+1)-W(0,n,m+1)\ge0
\end{align}
and
\begin{align}
\psi(n,k,m):=W(k-1,n,m+1)-W(1,n,m+1)\ge0.
\end{align}

By Ferroni \cite[Remark 3.9]{Fer}, $W(\ell,n,m)=W(n-m-\ell,n,m)$, then we have \[\varphi(n,k,m)=W(n-m-k+1,n,m+1)-W(0,n,m+1).\] Since $m<k\le\frac{n}{2}$, we find $n-(m+1)\ge n-m-k+1>0$. By the same arguments in the proof of \eqref{ww}, we can conclude that $W(n-m-k+1,n,m+1)\ge W(0,n,m+1)$. Thus $\varphi(n,k,m)\ge0$.

Similarly, since $n-(m+1)\ge k-1>1$, we can also utilize the same arguments in the proof of \eqref{ww} to show that $W(k-1,n,m+1)\ge W(1,n,m+1)$. That is, $\psi(n,k,m)\ge0$.
This completes the proof. \qed

\begin{lem} \label{stil}
For positive integers $m,n$, we have
\[
{n+1\brack m+1}=\sum_{j=0}^n{j\brack m}\frac{n!}{j!}.
\]
\end{lem}

\pf Recall that ${n+1\brack m+1}$ represents the number of permutations on $[n+1]$ with $m+1$ cycles. Alternatively, we can first choose $j$ numbers from $1,2,\ldots,n$ to form $m$ cycles, there are ${n\choose j}{j\brack m}$ such ways, then the left $n-j$ numbers and the number $n+1$ form another cycle, there are $(n-j)!$ such ways.
\qed

Now we are ready to give a proof of Theorem \ref{spp}.

\noindent{\it Proof of Theorem \ref{spp}}.  By Lemma \ref{umk} and Lemma \ref{uk2}, it suffices to show that \[i(U_{2,n},t)-i(T_{k,n},t-1)-i(T_{k,n},t)\] has positive coefficients. By Lemma \ref{umk} again, it is enough to only consider $k\le \frac{n}{2}$.

We first simplify $[t^m]i(U_{2,n},t)$, $[t^m]i(T_{k,n},t)$ and $[t^m]i(T_{k,n},t-1)$  for $0\le m\le n$, separately.

By Ferroni \cite[Corollary 3.13]{Fer},
\[
W(\ell,n,m)=\sum_{j=0}^\ell\sum_{i=0}^{n-m}(-1)^{i+j}{n\choose j}{m+\ell-j-1\choose m-1}{j\brack j-i}{n-j\brack m+i-j}.
\]
Then
\begin{align}\label{w1}
W(1,n,m+1)=(m+1){n\brack m+1}-n{n-1\brack m}.
\end{align}
Since $A(m,1)=2^m-m-1$ and $A(m,0)=1$, by \eqref{ukn} and \eqref{w1}, we have
\begin{align}
&(n-1)![t^m]i(U_{2,n},t)\nonumber\\
&=W(0,n,m+1)A(m,1)
+W(1,n,m+1)A(m,0)\nonumber\\
&=2^m{n\brack m+1}-n{n-1\brack m}\label{mu2k}\\
&=2^m\left((n-1){n-1\brack m+1}+{n-1\brack m}\right)
-n\left((n-2){n-2\brack m}+{n-2\brack m-1}\right)\nonumber\\
&=(n-2)\left(2^m{n-1\brack m+1}-(n-1){n-2\brack m}\right)+2\left(2^{m-1}{n-1\brack m}-(n-1){n-2\brack m-1}\right)\nonumber\\
&\quad +2^m{n-1\brack m+1}-(n-2){n-2\brack m}+(n-2){n-2\brack m-1}\nonumber\\
&\ge(n-2)\left(2^m{n-1\brack m+1}-(n-1){n-2\brack m}\right)+2\left(2^{m-1}{n-1\brack m}-(n-1){n-2\brack m-1}\right)\nonumber\\
&\quad+2^m{n-1\brack m+1}-(n-2){n-2\brack m}.\label{newu2k}
\end{align}

Replacing $t$ with $t-1$ in  \eqref{ehr} and extracting  the coefficient of $t^m$, we find
\begin{align}
&(n-1)!\cdot[t^m]i(T_{k,n},t-1)\nonumber\\
&= \sum_{j=0}^{k-1}\sum_{l=0}^{j}  \frac{(k-1)!}{j!}\binom{n-k-1+j}{j}  {j\brack l} {n-k\brack m-l}\label{eq53}\\
&= \sum_{j=0}^{k-1}\sum_{l=0}^{j}  \frac{(k-1)!}{j!}\frac{n-k-1+j}{n-k-1}\binom{n-k-2+j}{j}  {j\brack l} {n-k\brack m-l}.\nonumber
\end{align}
Since $n-k-1+j\le n-2$ and $\frac{n-k-1+j}{n-k-1}\le\frac{n-2}{n-k-1}\le 2$, we have
\begin{align*}
\frac{n-k-1+j}{n-k-1}{n-k\brack m-l}&=\frac{n-k-1+j}{n-k-1}\left((n-k-1){n-k-1\brack m-l}+{n-k-1\brack m-l-1}\right)\\
&\le (n-2){n-k-1\brack m-l}+2{n-k-1\brack m-l-1}.
\end{align*}
Then we conclude that
\begin{align}
&(n-1)!\cdot[t^m]i(T_{k,n},t-1)\nonumber\\
&\le \sum_{j=0}^{k-1}\sum_{l=0}^{j}  \frac{(k-1)!}{j!}\binom{n-k-2+j}{j}\left((n-2){n-k-1\brack m-l}+2{n-k-1\brack m-l-1}\right)  {j\brack l}.\label{newtm1}
\end{align}

Similarly, by  \eqref{ehr},
\begin{align}
&(n-1)!\cdot[t^m]i(T_{k,n},t)\nonumber\\
&= \sum_{j=0}^{k-1}\sum_{l=0}^{j}  \frac{(k-1)!}{j!}\binom{n-k-1+j}{j}  {j+1\brack l+1} {n-k+1\brack m-l+1}\label{weifangs}\\
&=\sum_{j=0}^{k-1}\sum_{l=0}^{j}  \frac{(k-1)!}{j!}\frac{n-k-1+j}{n-k-1}\binom{n-k-2+j}{j}  {j+1\brack l+1} {n-k+1\brack m-l+1}.\nonumber
\end{align}
Since $\frac{n-k-1+j}{n-k-1}(n-k)
=(n-2)+(j-k+2+\frac{j}{n-k-1})$ and $\frac{n-k-1+j}{n-k-1}\le 2$, we have
\begin{align*}
&\frac{n-k-1+j}{n-k-1}{n-k+1\brack m-l+1}\\
&=\frac{n-k-1+j}{n-k-1}\left((n-k){n-k\brack m-l+1}+{n-k\brack m-l}\right)\\
&\le (n-2){n-k\brack m-l+1}+\left(j-k+2+\frac{j}{n-k-1}\right){n-k\brack m-l+1}+2{n-k\brack m-l}.
\end{align*}
Therefore,  we find
\begin{align}
&(n-1)!\cdot[t^m]i(T_{k,n},t)\nonumber\\
&\le \sum_{j=0}^{k-1}\sum_{l=0}^{j}  \frac{(k-1)!}{j!}\binom{n-k-2+j}{j}\left((n-2){n-k\brack m-l+1}+2{n-k\brack m-l}\right){j+1\brack l+1}+h(n,m),\label{hnmm}
\end{align}
where
\[
h(n,m)=\sum_{j=0}^{k-1}\sum_{l=0}^{j}  \frac{(k-1)!}{j!}\binom{n-k-2+j}{j}\left(j-k+2+\frac{j}{n-k-1}\right){n-k\brack m-l+1}{j+1\brack l+1}.
\]

Consequently,  by \eqref{mu2k}, \eqref{eq53} and \eqref{weifangs}, we have
\begin{align*}
&[t^m](i(U_{2,n},t)-i(T_{k,n},t-1)-i(T_{k,n},t))\\
&= \frac{1}{(n-1)!}\left(2^m{n\brack m+1}-n{n-1\brack m}\right)\\
&\quad-\frac{1}{(n-1)!} \sum_{j=0}^{k-1}\sum_{l=0}^{j}  \frac{(k-1)!}{j!}\binom{n-k-1+j}{j}  \left({j\brack l} {n-k\brack m-l}+{j+1\brack l+1} {n-k+1\brack m-l+1}\right).
\end{align*}
Denote by
\begin{align*}
f(n,m)&=2^m{n\brack m+1}-n{n-1\brack m}\\
&\quad- \sum_{j=0}^{k-1}\sum_{l=0}^{j}  \frac{(k-1)!}{j!}\binom{n-k-1+j}{j} \left( {j\brack l} {n-k\brack m-l}+ {j+1\brack l+1} {n-k+1\brack m-l+1}\right).
\end{align*}
We aim to show that $f(n,m)\ge0$ for  $n\ge4$ and $m\ge0$ by induction on $n$.

It is easy to check that $f(4,m)\ge0$ for any $m\ge0$. Moreover, since
\[
\sum_{j=0}^{k-1}
\binom{n-k-1+j}{j}=\frac{k}{n-k}{n-1\choose k},
\]
we have
\begin{align*}
f(n,0)&=(n-1)!-(k-1)!(n-k)!-(n-k)!(k-1)!\sum_{j=1}^{k-1}
\binom{n-k-1+j}{j}\\
&=(n-1)!-(k-1)!(n-k)!-(n-k)!(k-1)!\left(\frac{k}{n-k}{n-1\choose k}-1\right)\\
&=0.
\end{align*}

Assume that $f(n-1,m)\ge0$ for $n\ge5$ and  $m\ge0$. For $m\ge1$, by \eqref{newu2k}, \eqref{newtm1} and \eqref{hnmm}, we derive that
\begin{align*}
f(n,m)
&\ge(n-2)\left(2^m{n-1\brack m+1}-(n-1){n-2\brack m}\right)+2\left(2^{m-1}{n-1\brack m}-(n-1){n-2\brack m-1}\right)\\
&\quad-\sum_{j=0}^{k-1}\sum_{l=0}^{j}  \frac{(k-1)!}{j!}\binom{n-k-2+j}{j}\left((n-2){n-k-1\brack m-l}+2{n-k-1\brack m-l-1}\right)  {j\brack l}\\
&\quad-\sum_{j=0}^{k-1}\sum_{l=0}^{j}  \frac{(k-1)!}{j!}\binom{n-k-2+j}{j}\left((n-2){n-k\brack m-l+1}+2{n-k\brack m-l}\right){j+1\brack l+1}\\
&\quad+2^m{n-1\brack m+1}-(n-2){n-2\brack m}-h(n,m)\\
&=(n-2)f(n-1,m)+2f(n-1,m-1)+2^m{n-1\brack m+1}-(n-2){n-2\brack m}-h(n,m)\\
&\ge 2^m{n-1\brack m+1}-(n-2){n-2\brack m}-h(n,m).
\end{align*}

To complete the proof, let
\[
g(n,m)=2^m{n-1\brack m+1}-(n-2){n-2\brack m}-h(n,m).
\]
We aim to prove that $g(n,m)\ge0$ for any $n\ge 4$ and $m\ge0$ by induction on $n$.

It's easy to check $g(4,m)= 0$ for any $m\ge0$ and $g(n,-1)=0$. Assume that $g(n-1,m)\ge0$ for $n\ge5$ and $m\ge0$. Then we have
\begin{align*}
g(n,m)&=2^m {n-1\brack m+1}-(n-2){n-2\brack m}\\
&\quad-\sum_{j=0}^{k-1} \sum_{l=0}^{j} \frac{(k-1)!}{j!}\left(j-k+2+\frac{j}{n-k-1}\right)\binom{n-k-2+j}{j}  {j+1\brack l+1} {n-k\brack m-l+1}\\
&=2^m\left( (n-2){n-2\brack m+1}+{n-2\brack m}\right)-(n-2)\left( (n-3){n-3\brack m}+{n-3\brack m-1}\right)\\
&\quad-\sum_{j=0}^{k-1} \sum_{l=0}^{j} \frac{(k-1)!}{j!}\left(j-k+2+\frac{j}{n-k-1}\right)\\
&\quad\times\frac{n-k-2+j}{n-k-2}\binom{n-k-3+j}{j}   \left((n-k-1){n-k-1\brack m-l+1}+{n-k-1\brack m-l}\right){j+1\brack l+1}\\
&\ge (n-2)\left(2^m{n-2\brack m+1}-(n-3){n-3\brack m}\right)+2\left(2^{m-1}{n-2\brack m}-(n-3){n-3\brack m-1}\right)\\
&\quad-\sum_{j=0}^{k-1} \sum_{l=0}^{j} \frac{(k-1)!}{j!}\left(j-k+2+\frac{j}{n-k-2}\right)\\
&\quad\times\binom{n-k-3+j}{j} \left((n-2){n-k-1\brack m-l+1}+2{n-k-1\brack m-l}\right) {j+1\brack l+1} \\
&\ge (n-2)g(n-1,m)+2g(n-1,m-1)\\
&\ge 0,
\end{align*}
where   the third step follows from the relations \[\left(j-k+2+\frac{j}{n-k-1}\right)
\frac{n-k-2+j}{n-k-2}(n-k-1)
\le \left(j-k+2+\frac{j}{n-k-2}\right)(n-2)\]
and
\[
\left(j-k+2+\frac{j}{n-k-1}\right)\frac{n-k-2+j}{n-k-2}
\le2\left(j-k+2+\frac{j}{n-k-2}\right).
\]
This completes the proof. \qed

\vspace{.2cm}
\noindent{\bf Acknowledgments.}
The authors wish to thank the referees for their valuable comments and
suggestions.
 We are grateful to Shaoshi Chen,   Peter Guo, Lisa Sun,  Matthew Xie, Sherry Yan and Arthur Yang for helpful conversations. Yao Li would like to thank the Research Experience for Undergraduates (REU) program ``Sparklet'' of the Math Department at Sichuan University. This work was supported by  the National Natural  Science Foundation of China  (Grant No. 11971250, 12071320) and Sichuan Science and Technology Program (Grant No. 2020YJ0006).

\footnotesize

D{\scriptsize EPARTMENT OF} M{\scriptsize ATHEMATICS}, S{\scriptsize ICHUAN} U{\scriptsize NIVERSITY}, C{\scriptsize HENGDU} 610064, P.R. C{\scriptsize HINA.}

Email address:
fan@scu.edu.cn (N.J.Y. Fan)

Email address:
liyaao@stu.scu.edu.cn (Y. Li)

\end{document}